\documentclass[12pt]{amsart}

\usepackage{amsmath,amssymb, amscd, stmaryrd}
\usepackage{mathtools}
\usepackage{fullpage}
\usepackage{enumerate}
\usepackage{framed}
\usepackage[dvips,dvipdf]{graphicx}
\usepackage[usenames,dvipsnames]{color}
\usepackage{tabularx}
\usepackage{array}
\usepackage{multirow}

\usepackage{amsfonts}
\usepackage{epstopdf}
\usepackage{subfigure,caption,float}
\usepackage{enumerate}
\usepackage{cleveref}
\usepackage{tabularx}

\newtheorem{theorem}{Theorem}[section]

\theoremstyle{definition}

\theoremstyle{definition}

\theoremstyle{remark}

\newtheorem{exmp}[theorem]{Example}

\title{Computation of highly oscillatory integrals \\ using a Fourier extension approximation}

\author[A. Anand]{Akash Anand} 
\address{Akash Anand, Department of
	Mathematics and Statistics, Indian Institute of Technology, Kanpur, UP 208016}
\email{akasha@iitk.ac.in}

\author[D. Dhiman]{Damini Dhiman} 
\address{Damini Dhiman, Department of
	Mathematics and Statistics, Indian Institute of Technology, Kanpur, UP 208016}
\email{damini@iitk.ac.in}


\begin{document}

\begin{abstract} 
The numerical evaluation of integrals of the form
\begin{align*}
	\int_a^b f(x) e^{ikg(x)}\,dx
\end{align*}
is an important problem in scientific computing with significant applications in many branches of applied mathematics, science and engineering. The  numerical approximation of such integrals using classical quadratures can be prohibitively expensive at high oscillation frequency ($k \gg 1$) as the number of quadrature points needed for achieving a reasonable accuracy must grow proportionally to $k$. 
To address this significant computational challenge, starting with Filon in 1930, several specialized quadratures have been developed to compute such oscillatory integrals efficiently. A crucial element in such Filon-type quadrature is the accurate evaluation of certain moments which poses a significant challenge when non-linear phase functions $g$ are involved.
		
In this paper, we propose an equispaced-grid Filon-type quadrature for computing such highly oscillatory integrals that utilizes a Fourier extension of the slowly varying envelope $f$. This strategy is primarily aimed at significantly simplifying the moment calculations, even when the phase function has stationary points. 
Moreover, the proposed approach can also handle certain integrable singularities in the integrand. 
We analyze the scheme to theoretically establish high-order convergence rates. We also include a wide variety of numerical experiments, including oscillatory integrals with algebraic and logarithmic singularities, to demonstrate the performance of the quadrature.

\end{abstract}

\maketitle

\section{Introduction}
The evaluation of highly oscillatory integrals plays a significant role in many branches of applied mathematics. The accuracy of their numerical approximation using a classical quadrature suffers significantly when the oscillation frequency is much more than the number of quadrature points. Since the initial attempt by Filon \cite{filon1930iii} to solve this difficulty, numerous approaches have been developed \cite{luke1954computation, flinn1960modification, tuck1967simple, iserles2004numerical, iserles2004quadrature}. The Filon quadrature for approximating the integral of the form 
\begin{align} \label{genprob}
\int_a^b f(x)e^{i k g(x)}\,dx,
\end{align}
relies on the interpolating polynomial $p$ satisfying $p(x_j)=f(x_j)$ for interpolation nodes $x_0, x_1, \ldots, x_{n},$ to obtain the approximation
\begin{align} \label{filonapprox}
\int_a^b f(x)e^{i k g(x)}\,dx  \approx  \int_a^b p(x)e^{i k g(x)}\,dx.
\end{align}
Both the Filon method and its generalisation  \cite{luke1954computation, flinn1960modification, tuck1967simple} where $f$ is approximated by splines, are shown to be efficient for suitably smooth functions provided the moments $\int_a^b x^j e^{ik g(x)}\,dx$ can be explicitly calculated. While obtaining these moments is easy for linear oscillators, that is, $g(x) = \alpha x + \beta,\ \alpha, \beta\in\mathbb{R}$, in general, the accurate moment calculations pose a significant challenge. For example, The Filon-Clenshaw-Curtis quadrature \cite{dominguez2011stability}, that uses $\cos(j\pi/N), j = 0,1,\dots, N$ as nodes and Chebyshev approximations for $f$, obtain the underlying moments recursively using multiple recurrence relations. The Filon-Clenshaw-Curtis and its variations \cite{dominguez2013filon, dominguez2014filon, li2019efficient, kang2011efficient} have been developed to efficiently handle integrable algebraic  and logarithmic  singularities in the integrand.

Several Filon-like approaches use the change of variable $y=g(x)$ to transform the integral to the one with linear oscillations
\begin{align}\label{integral}
\int_a^b f(x)e^{i k g(x)}\,dx = \int_{g(a)}^{g(b)} \frac{1}{g'(g^{-1}(y))}f(g^{-1}(y))e^{i k y}\,dy
\end{align}
under the assumption that $g'(x) > 0$ for all $x \in (a,b)$ or $g'(x)<0$ for all $x\in(a,b)$. Indeed, if $g'(x)$ has finitely many zeros in $(a,b)$, say $x_j, j = 1,2,\ldots,J$, then, by setting $x_0=a$ and $x_{J+1}=b$, we have 
\begin{align*}
\int_a^b f(x)e^{i k g(x)}\,dx =  \sum_{j=0}^{J} \int_{x_j}^{x_{j+1}} f(x)e^{i k g(x)}\,dx
\end{align*}
where, for each $j$, the integral over $[x_j,x_{j+1}]$ can be transformed into the one with linear oscillations.

In this paper, we propose a Filon-type quadrature for oscillatory integrals with equispaced quadrature nodes. We show that this approach not only handles algebraic and logarithmic singularities in the integrand but also significantly simplifies the solution of the underlying moment problem. In fact, for several special singular weights, the proposed approach yields moment problems that can be solved analytically. Moreover, our approach requires only discrete functional data on a uniform grid.  The fact that approximating the integral to high order does not require explicit knowledge of the derivatives is a significant strength of this approach.

The method to yield a more tractable moment problem relies on the interpolating trigonometric polynomial approximation of the smooth non-oscillatory part of the integrand. However, in general, such an approximation is inflicted with
Gibb's oscillations that manifest due to the discontinuous periodic extension of the function being approximated. To overcome this difficulty, we first construct a periodic extension of the function, which is then approximated by the interpolating trigonometric polynomial using FFT. We discuss the periodic extension employed in the paper in \cref{PeriodicExt}. In \cref{MomentGen}, we delve into the moment problem and discuss the solution of several important instances. 
 In \cref{ErrorAnalysis}, we present a convergence analysis. A variety of numerical experiments are presented in \cref{NumericalResults} to demonstrate the performance of the proposed quadrature. 

\section{Methodology}

Motivated by the need to compute the integrals such as those that appear on the right-hand side of the \cref{integral}, we consider the problem of approximation of the integrals
\begin{align} \label{WeightGenProb}
I^{[a,b]}_k (w,f) = \int_a^b w(x) f(x)e^{i k x}\,dx,
\end{align}
where $w \in L^{1}(a,b)$ is the weight function and $f$ is a smooth function on $[a,b]$. 

\subsection{The basic idea} \label{Idea}

For clarity, we introduce the main ideas assuming that the $f$ is known everywhere in the interval $[a,b]$ and its derivatives are available at $x=a$ and $x=b$. At the heart of this procedure is the step where we approximate $f$ by a trigonometric polynomial 
\begin{align*}
\sum_{\ell=-n}^{n-1}c_\ell(f) e^{2\pi i\ell (x-a)/\xi}
\end{align*}
of period $\xi \ge b-a$. The primary motivation for doing so is to simplify the underlying moment calculations. Indeed, we see that 
\begin{align} \label{ContinuousApprox}
I^{[a,b]}_k (w,f) \approx \int_a^b w(x) \left(\sum_{\ell=-n}^{n-1}c_\ell(f) e^{2\pi i\ell (x-a)/\xi}\right) e^{i k x}\,dx = e^{ika} \sum_{\ell=-n}^{n-1} c_{\ell}(f) W^{[a,b]}_{k,\ell}(w) 
\end{align}
where
\begin{align} \label{Moments}
c_{\ell}(f) = \frac{1}{\xi}\int_a^{a+\xi} f(x)e^{-2\pi i \ell(x-a)/\xi}\,dx \quad \text{and}\quad W^{[a,b]}_{k,\ell}(w) = \int_a^bw(x)e^{i(k\xi+2\pi \ell) (x-a)/\xi}\,dx.
\end{align}
For many important weights $w$, the moments $W^{[a,b]}_{k,\ell}(w)$ can be obtained analytically as we see later in \cref{MomentGen}.

The natural choice for the approximating trigonometric polynomial of $f$ is to use its truncated Fourier series. However, if we use the period $\xi=b-a$, it fails to converge uniformly if the $(b-a)$-periodic extension of $f$ is discontinuous, that is, if $f(a)\ne f(b)$, due to Gibb's phenomenon. 
Instead of directly approximating $f$, we utilize the Fourier series of a periodic function with a larger period that coincides with $f$ on the interval $[a,b]$ to eliminate Gibb's oscillations. 
In the next sub-section, we discuss one such construction that we use for numerical experiments in this paper.

\subsection{A periodic extension:} \label{PeriodicExt}

Among several possible periodic extensions for a given function, for simplicity of implementation and analysis, we utilize the construction based on the two-point Hermite polynomial  \cite{anand2019fourier, huybrechs2006evaluation}.
Given $f$ on $[a,b]$, we introduce a periodic function $f^e$ of period $\xi = 2(b-a)$ given by
\begin{align} \label{feGen}
f^e(x) = 
\begin{dcases}
f(x), & x \in [a,b], \\
p(x), & x \in [b,2b-a),
\end{dcases}
\end{align}
where the polynomial $p$ is chosen as
\begin{align} \label{Polynomialp}
p(x) = \sum_{m=0}^r f^{(m)}(b) p_{m,r}^{[b,2b-a]}(x) +  \sum_{m=0}^r f^{(m)}(a) p_{m,r}^{[2b-a,b]}(x)
\end{align}
with
\begin{align*}
p_{m,r}^{[t_1,t_2]}(t) &=  \frac{1}{m!} (t-t_1)^m \left(\frac{t-t_2}{t_1-t_2}\right)^{r+1}  \sum_{n=0}^{r-m} \binom{r+n}{r}\left(\frac{t-t_1}{t_2-t_1}\right)^n.
\end{align*}


\subsection{The Discrete Problem:} \label{DiscreteProb}
In many practical applications, only discrete functional data is available. In this regard, we consider the equispaced grid of size $n+1$ on the interval $[a,b]$ with grid points $x_j = a+(b-a)j/n$ where the functional data  $f_j, j = 0,1,\dots,n$ is assumed to be available. We introduce the grid function $f_{r,q}^e$ given as
\begin{align} \label{Discretefegen}
(f_{r,q}^e)_j = 
\begin{dcases}
f_j, & j=0,\dots,n, \\
p_{r,q}(a+(b-a)j/n), & j=n+1,\dots,2n,\\
\end{dcases}
\end{align}
where
\begin{align} \label{discreteExtn}
p_{r,q}(x) = \sum_{m=0}^r \mathcal{D}^{m,-}_{n,q}(f)(x_n) p_{m,r}^{[b,2b-a]}(x)+ \sum_{m=0}^r \mathcal{D}^{m,+}_{n,q}(f)(x_0) p_{m,r}^{[2b-a,b]}(x),
\end{align}
for $1\le m \le r$, the operators $\mathcal{D}^{m,+}_{n,q}$ and $\mathcal{D}^{m,-}_{n,q}$ respectively denote the forward and backward finite difference derivative operators for approximating $f^{(m)}$ of the form
\begin{align*}
\mathcal{D}^{m,\pm}_{n,q}(f)(x_j) = (\pm n)^m\left( \sum_{\ell = 0}^{m+q-1} (a_q^m)_\ell f_{j\pm \ell}\right) 
 \end{align*}
with appropriately chosen constants $(a_q^m)_\ell$ so that the order of accuracy is $q$. The discrete Fourier coefficients are obtained as
\begin{align} \label{DiscreteCoeff}
d_{\ell,n}(f_{r,q}^e) = \frac{1}{2n} \sum_{j=0}^{2n-1} (f_{r,q}^e)_j e^{- \pi i \ell j/n}
\end{align}
for $\ell = -n,\dots,n-1$.
Finally, the integration scheme for approximation of $I_k^{[a,b]}(w,f)$ is given by
\begin{align} \label{DiscreteApprox}
I_k^{[a,b]}(w,f)\approx I_{k,n}^{[a,b]}(w,f) := e^{ika}\sum_{\ell=-n}^{n-1}d_{\ell,n}(f_{r,q}^e)W_{k,\ell}^{[a,b]}(w).
\end{align}

\subsection{The moment problem:} \label{MomentGen}

The numerical integration scheme requires obtaining the moments given by
\begin{align}
W^{[a,b]}_{k,\ell}(w) &=
\int_a^bw(x)e^{\pi i\ell (x-a)/(b-a)}e^{ik(x-a)}\,dx 
= (b-a) \int_0^1 w(a+(b-a)u)e^{i((b-a)k+\pi\ell)u}\,du \nonumber\\
&= \frac{b-a}{2} e^{i((b-a)k+\pi\ell)/2} \int_{-1}^{1} w\left(\frac{a+b}{2}+\frac{b-a}{2}v \right) e^{i((b-a)k+\pi\ell)v/2}\,dv \nonumber
\end{align}
The following are solutions to some instances of the moment problem corresponding to important classes of weights of interest.
\begin{exmp}
For the simplest case where $w(x) =1$,
we have
\begin{align*}
W^{[a,b]}_{k,\ell}
&= \begin{cases}
b-a, & \ell = -(b-a)k/\pi, \\
\dfrac{2(b-a)}{(b-a)k+\pi\ell} e^{i((b-a)k+\pi\ell)/2} \sin(((b-a)k+\pi\ell)/2), & \ell \ne -(b-a)k/\pi.
\end{cases}
\end{align*}
\end{exmp}

\begin{exmp}
For one sided algebraic singularity where $w(x) = (x-a)^{\beta},$ $\beta \in (-1,0)$, the moments are given by
\begin{align*}
W^{[a,b]}_{k,\ell} 
&= \begin{cases}
\dfrac{(b-a)^{1+\beta}}{1+\beta}, & \ell=-(b-a)k/\pi, \\
 \left(\dfrac{i(b-a)}{(b-a)k+\pi\ell}\right)^{1+\beta} (\Gamma(1+\beta)-\Gamma(1+\beta, -i((b-a)k+\pi\ell))), & \ell\ne -(b-a)k/\pi.   
\end{cases}
\end{align*}
\end{exmp}

\begin{exmp}
For $w(x) = (x-a)^{\beta}(b-x)^{\beta}$ with $\beta \in (-1,0)$, the weak algebraic singularity at both end-points of the same order $\beta$, we have 
\begin{align*}
&W^{[a,b]}_{k,\ell}
= \sqrt{\pi}\ \Gamma(1+\beta)(b-a)^{1+2\beta} e^{i((b-a)k+\pi\ell)/2} ((b-a)k+\pi\ell)^{-1/2-\beta} J_{1/2+\beta}(((b-a)k+\pi\ell)/2)
\end{align*}
provided $\ell \ne -(b-a)k/\pi$, whereas 
\begin{align*}
W^{[a,b]}_{k,\ell}= \left(\frac{b-a}{2}\right)^{1+2\beta}\frac{\sqrt{\pi}\hspace{0.1cm} \Gamma(1+\beta)}{\Gamma\left(\frac{3}{2}+\beta\right)},\\
\end{align*}
if $\ell = -(b-a)k/\pi$. In particular, for $\beta = -1/2$, we get
\begin{align*}
W^{[a,b]}_{k,\ell} =\pi e^{i((b-a)k+\pi\ell)/2}  J_{0}(((b-a)k+\pi\ell)/2).
\end{align*}
\end{exmp}

\begin{exmp}
For $w(x) = (x-a)^\alpha(b-x)^\beta$, for $\alpha, \beta \in(-1,0)$, when the singularity is present at both end-points with different exponents, then
\begin{align*}
W^{[a,b]}_{k,\ell} 
&= \frac{b-a}{2} e^{i((b-a)k+\pi\ell)/2}\int_{-1}^{1}\left(\frac{b-a}{2}\right)^{\alpha+\beta}(1+v)^\alpha(1-v)^\beta e^{i((b-a)k+\pi\ell)v/2}\,dv \\
&= \frac{(b-a)^{\alpha+\beta+1}}{\Gamma(2+\alpha+\beta)}\Gamma(1+\alpha)\Gamma(1+\beta) _1F_1(1+\alpha,2+\alpha+\beta,i((b-a)k+\pi\ell)).
\end{align*}
\end{exmp}

\begin{exmp}
For $w(x) = \log(x-a)$, when the logarithmic singularity is present at the left end-point, then
\begin{align*}
W^{[a,b]}_{k,\ell}
=  \frac{i(b-a)}{(a-b)k-\pi\ell}(&\gamma+\Gamma(0,-i((b-a)k+\pi\ell))+(-1+e^{i((b-a)k+\pi\ell)})\log(b-a)\\
& +\log(-i((b-a)k+\pi\ell)),
\end{align*}
if $(b-a)k+\pi\ell \ne0$ and $W^{[a,b]}_{k,\ell} = (b-a)(\log(b-a)-1)$ otherwise. 
\end{exmp}

\section{Error Analysis} \label{ErrorAnalysis}

We begin by noting that, using the Fourier series of the periodic extension of $f$, that is, 
\begin{align*}
f(x)=
\sum_{\ell=-\infty}^{\infty}c_\ell(f^e)e^{\pi i\ell (x-a)/(b-a)},
\end{align*}
where
\begin{align*}
c_\ell(f^e)=\frac{1}{\xi}\int_{a}^{a+\xi}f^e(x)e^{-2\pi i\ell (x-a)/\xi}\,dx=\frac{1}{2(b-a)}\int_{a}^{2b-a}f^e(x)e^{-\pi i\ell (x-a)/(b-a)}\,dx,\\
\end{align*}
the exact integral can be expressed as
\begin{align*}
I^{[a,b]}_k(w,f) = e^{ika}
\sum_{\ell=-\infty}^{\infty}c_\ell(f^e)W^{[a,b]}_{k,\ell}(w).
\end{align*}
The error in the numerical integration, therefore, is given by 
\begin{align} 
&I^{[a,b]}_k(w,f) - I^{[a,b]}_{k,n}(w,f) = 
e^{ika} \left(\sum_{\ell=-\infty}^{\infty}c_{\ell}(f^e)W_{k,\ell}^{[a,b]}(w) -\sum_{\ell=-n}^{n-1}d_{\ell,n}(f_{r,q}^e)W_{k,\ell}^{[a,b]}(w) \right) = \nonumber\\
&e^{ika} \left(\sum_{\ell=-n}^{n-1} \sum_{\substack{m=-\infty \\ m \ne 0}}^\infty c_{\ell+2mn}(f^e)W_{k,\ell+2mn}^{[a,b]}(w)  +\sum_{\ell=-n}^{n-1}\left(c_{\ell}(f^e)-d_{\ell,n}(f_{r,q}^e)\right)W_{k,\ell}^{[a,b]}(w) \right) \label{eq:error_1}
\end{align}
Using
\begin{align*}
c_{\ell}(f^e)-d_{\ell,n}(f_{r,q}^e) 
= c_\ell(f^e-f_{r,q}^e) - \sum_{\substack{m=-\infty \\ m \ne 0}}^\infty c_{\ell+2mn}(f_{r,q}^e),
\end{align*}
we rewrite the second summation in \cref{eq:error_1} as
\begin{align} \label{eq:error_2}
    \sum_{\ell=-n}^{n-1} c_\ell(f^e-f_{r,q}^e)W_{k,\ell}^{[a,b]}(w) - \sum_{\ell=-n}^{n-1}\sum_{\substack{m=-\infty \\ m \ne 0}}^\infty c_{\ell+2mn}(f_{r,q}^e)W_{k,\ell}^{[a,b]}(w),
\end{align}
and, therefore, combining the expressions in \cref{eq:error_1} and \cref{eq:error_2}, the error reads
\begin{align} \label{eq:error_3}
    I^{[a,b]}_k(w,f) - I^{[a,b]}_{k,n}(w,f) = 
    e^{ika} \left(\sum_{\ell=-n}^{n-1} \sum_{\substack{m=-\infty \\ m \ne 0}}^\infty c_{\ell+2mn}(f^e)\left( W_{k,\ell+2mn}^{[a,b]}(w) - W_{k,\ell}^{[a,b]}(w) \right) + \right. &\nonumber\\ 
    \left.\sum_{\ell=-n}^{n-1} \sum_{\substack{m=-\infty \\ m \ne 0}}^\infty c_{\ell+2mn}(f^e-f_{r,q}^e) W_{k,\ell}^{[a,b]}(w) +\sum_{\ell=-n}^{n-1}  c_\ell(f^e-f_{r,q}^e)W_{k,\ell}^{[a,b]}(w) \right) &
\end{align}
To see that
\begin{align} \label{eq:error_4}
\sum_{\ell=-n}^{n-1}  c_\ell(f^e-f_{r,q}^e)W_{k,\ell}^{[a,b]}(w)  = -\sum_{\ell=-n}^{n-1} \sum_{\substack{m=-\infty \\ m \ne 0}}^\infty c_{\ell+2mn}(f^e-f_{r,q}^e) W_{k,\ell+2mn}^{[a,b]}(w)
\end{align}
we introduce
\begin{align*} 
    A(g)(z):=\int_a^{2b-a} w_0(y)g(z+y)e^{ik(y-a)}\,dy
\end{align*}
where $w_0$ is the  $2(b-a)$-periodic function with $w_0(x) = w(x)$ for $x\in [a,b]$ and $w_0(x) = 0$ for $x\in (b,2b-a)$. As
\begin{align*}
    &A(g)(z) = 
    \sum_{\ell=-\infty}^{\infty} c_\ell(g)e^{\pi i \ell z/(b-a)} W^{[a,b]}_{k,\ell}(w),
\end{align*}
 we have
\begin{align*} 
    &\sum_{\ell=-\infty}^{\infty} c_\ell(f^e-f_{r,q}^e) W^{[a,b]}_{k,\ell}(w) = A(f^e-f_{r,q}^e)(0) =0.
\end{align*}
Therefore, using \cref{eq:error_4}, and after simplification, the error in \cref{eq:error_3} is expressed as
\begin{align} \label{eq:error_5}
    &I^{[a,b]}_k(w,f) - I^{[a,b]}_{k,n}(w,f) = 
    e^{ika} \sum_{\ell=-n}^{n-1} \sum_{\substack{m=-\infty \\ m \ne 0}}^\infty c_{\ell+2mn}(f_{r,q}^e) \left( W_{k,\ell+2mn}^{[a,b]}(w) - W_{k,\ell}^{[a,b]}(w) \right).
\end{align}

\begin{figure}[t]
\centering
\subfigure[$w(x) = x^\beta,\ k = 10$ ]{\includegraphics[clip=true, trim=230 100 200 100, width=0.47\textwidth]{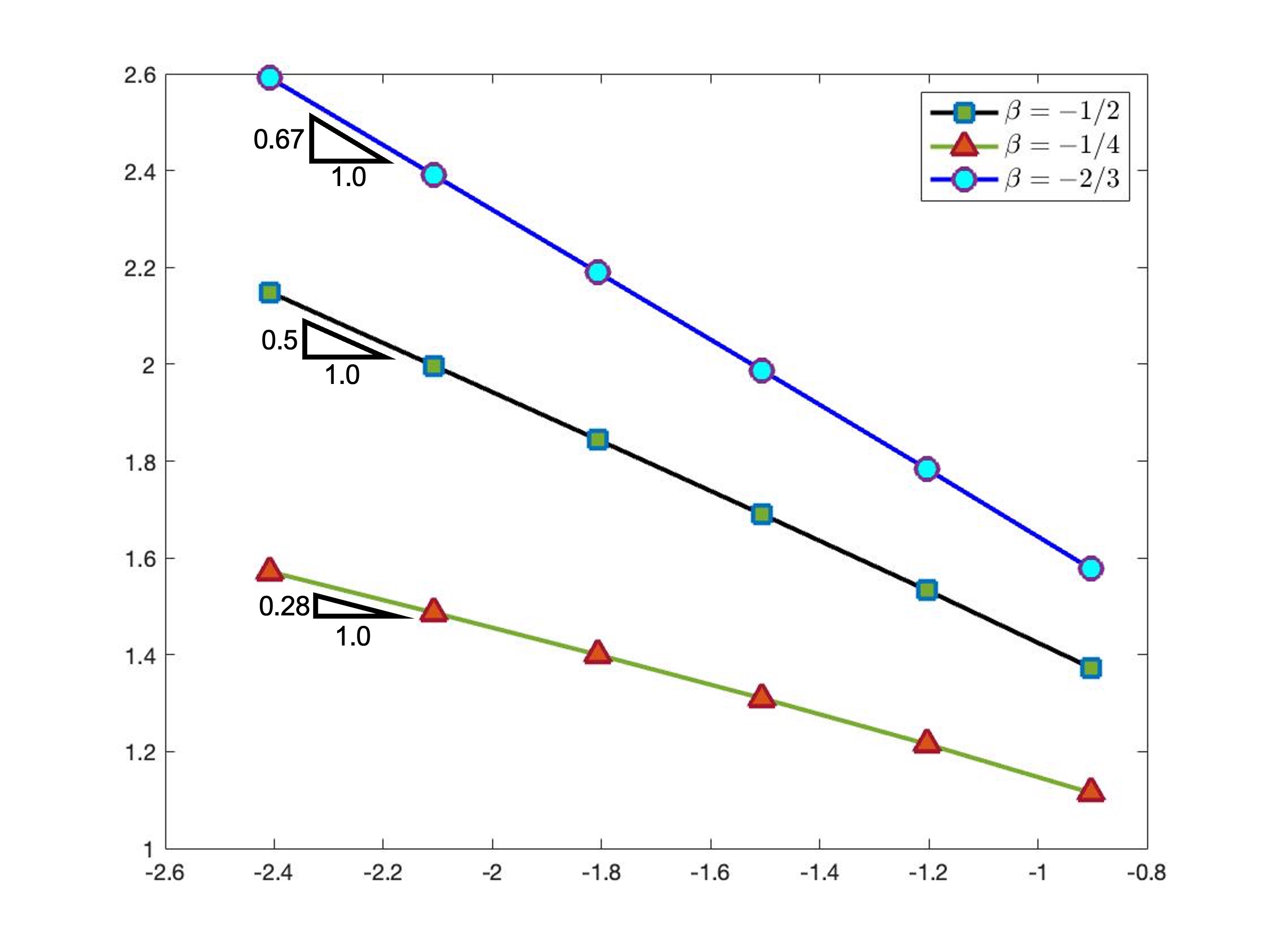} 
\label{fig:w1n}
} 
\subfigure[$w(x) = x^\alpha(1-x)^\beta, \ k = 10$]{\includegraphics[clip=true, trim=230 100 200 100, width=0.47\textwidth]{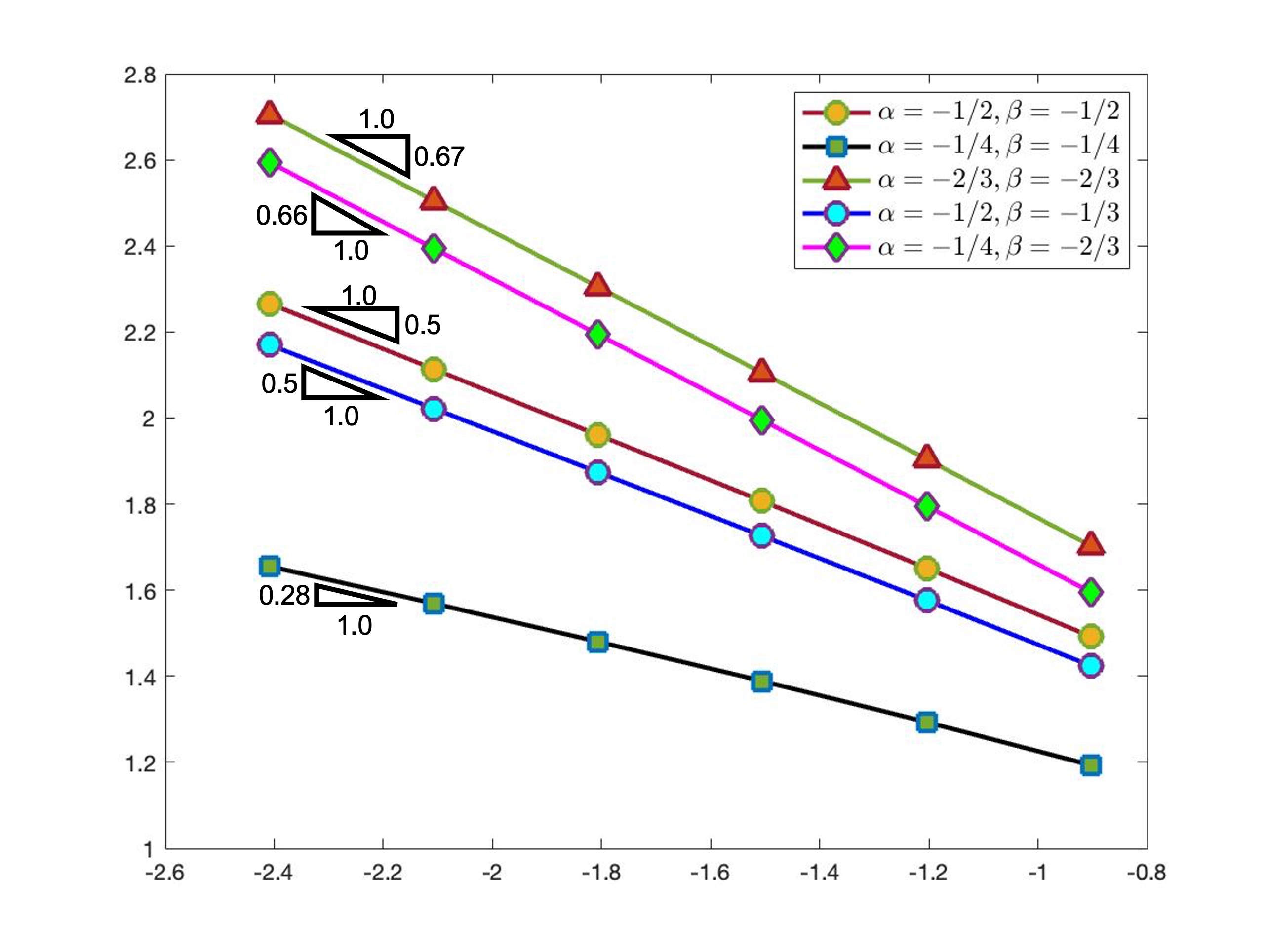} 
\label{fig:w2n}
} 
\subfigure[$w(x) = x^\beta, \ k = 100$ ]{\includegraphics[clip=true, trim=230 100 200 100, width=0.47\textwidth]{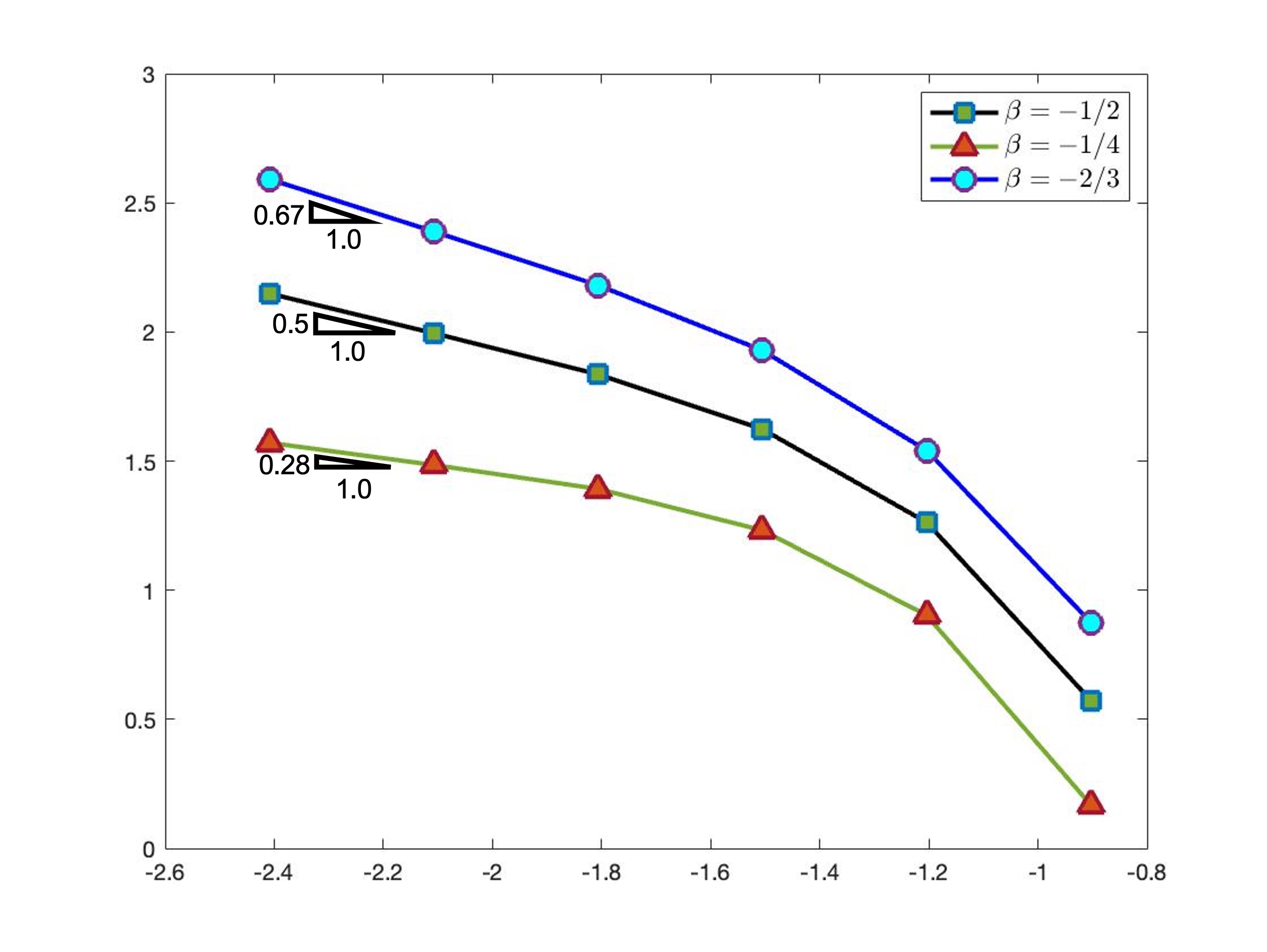} 
\label{fig:w11n}
} 
\subfigure[$w(x) = x^\alpha(1-x)^\beta, \ k = 100$]{\includegraphics[clip=true, trim=230 100 200 100, width=0.47\textwidth]{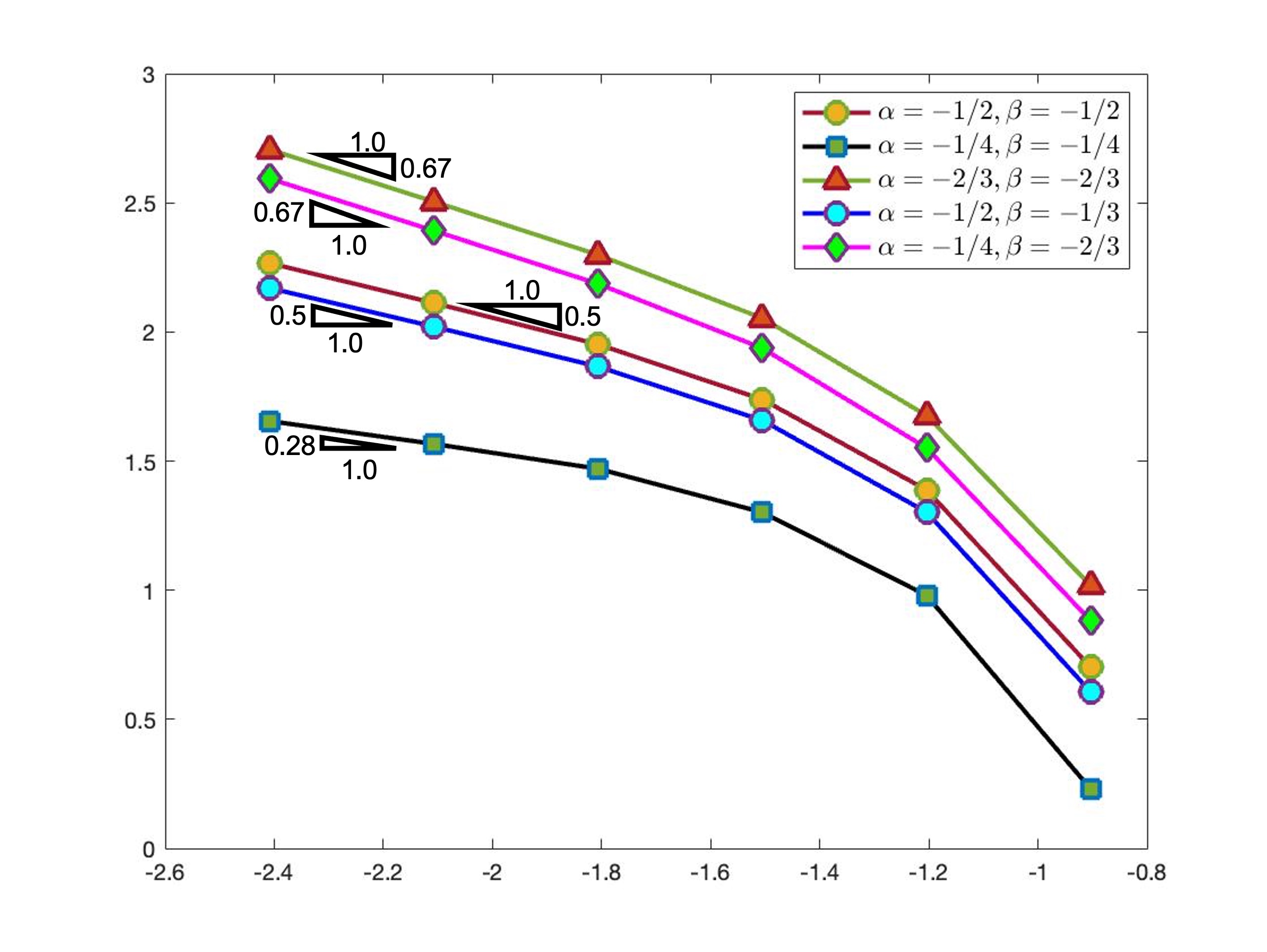} 
\label{fig:w22n}
} 
\caption{The plot of $\log_{10} W_{k,n,300}^{[a,b]}(w)$ against $\log_{10} \left(1/n\right)$}
\label{fig:wn}
\end{figure}

Thus, to obtain an error bound, we need an estimate on $\left|c_{\ell}(f_{r,q}^e)\right|$. In this direction, it is important to note that $f_{r,q}^e(x) = f^e(x) + n^{-q}g(x)$ where $g(x)=0$ for $x\in [a,b]$, $g\in C^2(b,2b-a)$ and $g(b)=g(2b-a)=0$. Thus, for $\ell\ne 0$, we have
\begin{align} \label{bound1}
   &c_{\ell}(f_{r,q}^e) 
   = c_{\ell}(f^e) +n^{-q}\frac{b-a}{2(\pi i \ell)^2} \left((-1)^\ell g'(b)-g'(2b-a) +\int_{b}^{2b-a}g''(x)e^{-\pi i\ell (x-a)/(b-a)}\,dx\right) 
\end{align}
and
\begin{align} \label{bound2}
   &c_{\ell}(f^e) =
   \frac{(b-a)^{r+1}}{2(\pi i \ell)^{r+2}} \Bigg((-1)^\ell \left(p^{(r+1)}(b)-f^{(r+1)}(b)\right) +f^{(r+1)}(a)-p^{(r+1)}(2b-a) \\
   &\quad\quad\quad\quad\quad +\int_{a}^{b}f^{(r+2)}(x)e^{-\pi i\ell (x-a)/(b-a)}\,dx+\int_{b}^{2b-a}p^{(r+2)}(x)e^{-\pi i\ell (x-a)/(b-a)}\,dx \Bigg). \nonumber
\end{align}
Therefore, using \cref{bound1} and \cref{bound2}, we conclude that there exits a positive constant $C$ such that
\begin{align}  \label{bound}
\left|c_{\ell}(f_{r,q}^e)\right| \le C \left( |\ell |^{-(r+2)} + n^{-q}  |\ell |^{-2} \right), \quad \ell\in \mathbb{Z}\setminus\{0\}.
\end{align}
Finally, combining \cref{eq:error_5} and \cref{bound}, we have
\begin{align} \label{errorbound}
&\left| I^{[a,b]}_k(w,f) - I^{[a,b]}_{k,n}(w,f)\right| \le 
C n^{-(\min\{r,q\}+2)} \sum_{\ell=-n}^{n-1}  \sum_{m=-\infty}^\infty \frac{\left| W_{k,\ell+2mn}^{[a,b]}(w) - W_{k,\ell}^{[a,b]}(w) \right|}{\left(2|m|-1\right)^2}.
\end{align}
We, therefore, see that the convergence rate for the approximations does depend on the weight function $w$. In this context, using the notation
\begin{align*}
    \mathcal{W}^{[a,b]}_{k,n,M}(w) := \sum_{\ell=-n}^{n-1}  \sum_{m=-M}^M \frac{\left| W_{k,\ell+2mn}^{[a,b]}(w) - W_{k,\ell}^{[a,b]}(w) \right|}{\left(2|m|-1\right)^2}
\end{align*}
we introduce a class of weight functions defined as
\begin{align*}
    \mathbb{W}_\gamma(a,b) := \left\{ w \in L^1(a,b) :\exists\  n_0 \in \mathbb{N},\ C>0, \mathcal{W}^{[a,b]}_{k,n,\infty}(w) \le Cn^\gamma,  \text{ for all }n \ge n_0 \right\}
\end{align*}
for which the approximations converge at the same rate for a given $f$. Obviously, $\mathbb{W}_{\gamma_1}(a,b) \subseteq \mathbb{W}_{\gamma_2}(a,b)$ for $\gamma_1 \le \gamma_2$. The numerical calculations plotted in \cref{fig:wn} indicate that the the weight $w = (x-a)^\beta, \beta\in (-1,0)$ belongs to $\mathbb{W}_{-\beta}(a,b)$ whereas $w = (x-a)^\alpha(b-x)^\beta, \alpha,\beta\in(-1,0)$ belongs to $\mathbb{W}_{\gamma}(a,b)$ with $\gamma = \max\{-\alpha,-\beta\}$.

The following theorem captures the essence of the error analysis of the proposed integration scheme.

\begin{theorem}
Let $w \in \mathbb{W}_\gamma(a,b)$ and $f \in C^{\infty}([a,b])$. Let  
\begin{align*}
I_{k,n}^{[a,b]}(w,f) := e^{ika}\sum_{\ell=-n}^{n-1}d_{\ell,n}(f_{r,q}^e)W_{k,\ell}^{[0,1]}(w), \quad q,n \in \mathbb{N}, \ \ r \in \mathbb{N}\cup\{0\},
\end{align*} 
where
\begin{align*}
    &f_{r,q}^e(x) =
    \begin{dcases}
        f(x), & x \in [a,b], \\
        \sum_{m=0}^r \mathcal{D}^{m,-}_{n,q}(f)(b) p_{m,r}^{[b,2b-a]}(x)+ \sum_{m=0}^r \mathcal{D}^{m,+}_{n,q}(f)(a) p_{m,r}^{[2b-a,b]}(x), & x \in (b,2b-a),
    \end{dcases} \\
&d_{\ell,n}(f_{r,q}^e) = \frac{1}{2n} \sum_{j=0}^{2n-1} f_{r,q}^e(a+(b-a)j/n) e^{- \pi i \ell j/n}, \ W^{[a,b]}_{k,\ell}(w) = \int_a^bw(x)e^{i(k(b-a)+\pi \ell) (x-a)/(b-a)}\,dx,
\end{align*}
be the approximation to the oscillatory integral $I_{k}^{[a,b]}(w,f)$.
Then, there exist $n_0 \in \mathbb{N}$ and a constant $C>0$ such that
\begin{align*}
    \left| I_k^{[a,b]}(w,f) - I_{k,n}^{[a,b]}(w,f)\right| \le Cn^{-(\min\{r,q\}+2-\gamma)},
\end{align*}
for all $n \ge n_0$.
\end{theorem}


The optimal choice $q=r$ is, therefore, recommended in practice for approximating the derivatives that results in order $r+2-\gamma$ rate of convergence for integral approximations.

\section{Numerical Results} \label{NumericalResults}

 In this section, we present a variety of numerical experiments to validate the performance of the proposed method. Unless otherwise stated, we use $q=r$ in quadrature approximations.  


\begin{exmp} \label{Eg1} To begin with, we test our method on an example used by David Levin in \cite{levin1982procedures} and consider the integral 
\begin{align*}
I_k = \int_{0}^{1}\sin(t)e^{ik(t+t^2)}\,dt.
\end{align*}
The change of variable $x = t+t^2$ yields $I_k = \int_0^2 f(x)e^{ikx}\,dx$ where
\begin{align*}
f(x) =  \sin\left(\frac{-1+\sqrt{4x+1}}{2}\right)\frac{1}{\sqrt{4x+1}}
\end{align*} 
is smooth on the interval $[0,2]$. Thus, the quadrature with weight function $w(x) = 1$ produces rapidly convergent approximations $I_{k,n}$ of $I_k$ where the convergence rate can be raised arbitrarily by increasing the parameter $r$. We studied the convergence at $k=100$, $k=500$ and $k = 1000$ for $r=1, 2, 3, 4,$ and presented the results in \cref{fig:exmp1} where we plot $\log_{10}|I_{k,n}-I_k|$ against $\log_{10}(1/ n)$ to gauge the speed of convergence through the slope of respective lines. The plots clearly show that approximations converge with the rate $r+2$.

\begin{figure}[t]
\begin{center}
\subfigure[$k=100$ ]{\includegraphics[clip=true, trim=240 130 240 140, width=0.31\textwidth]{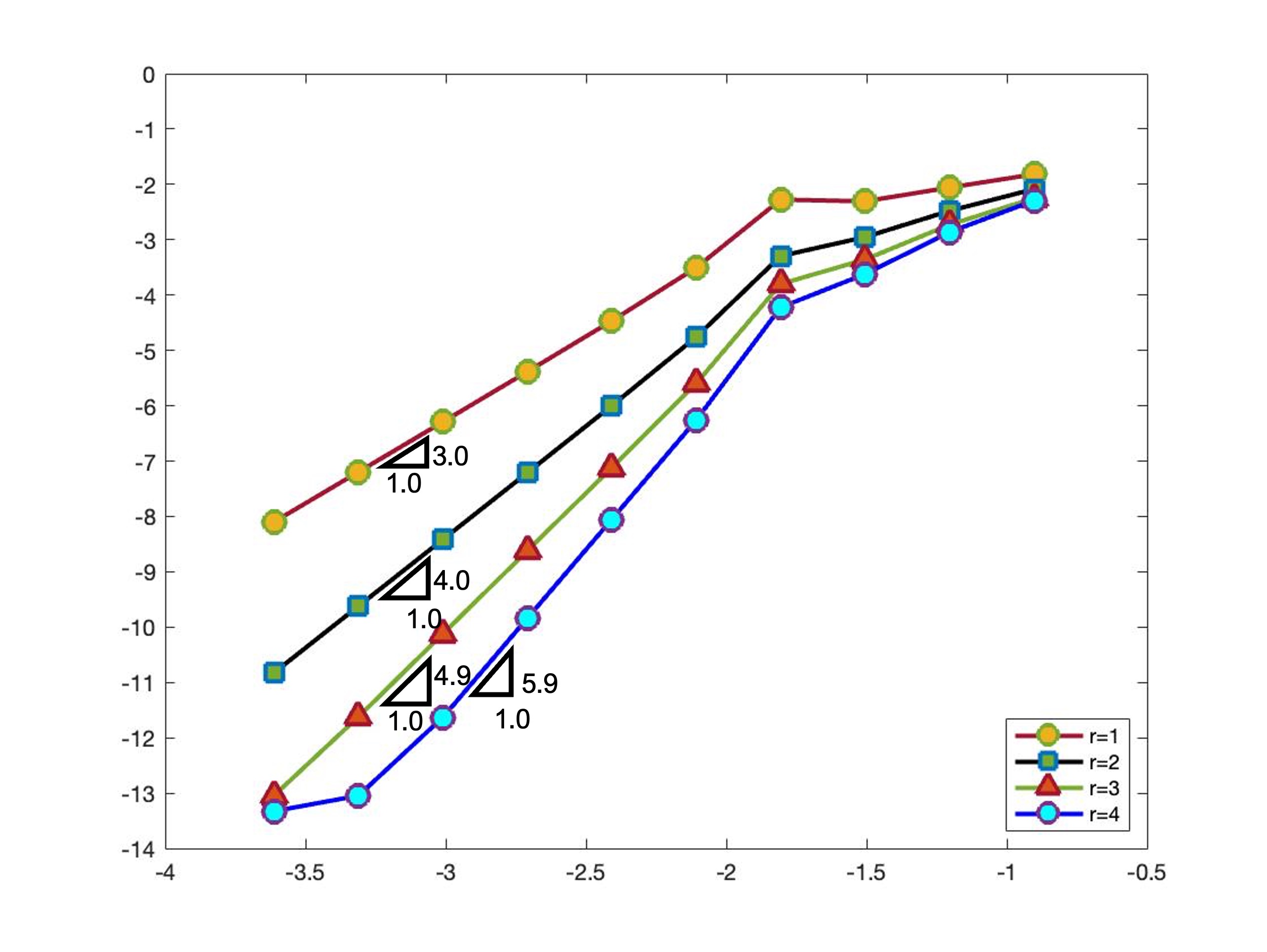} 
\label{fig:exmp1a}
} 
\subfigure[$k=500$]{\includegraphics[clip=true, trim=300 150 240 160, width=0.31\textwidth]{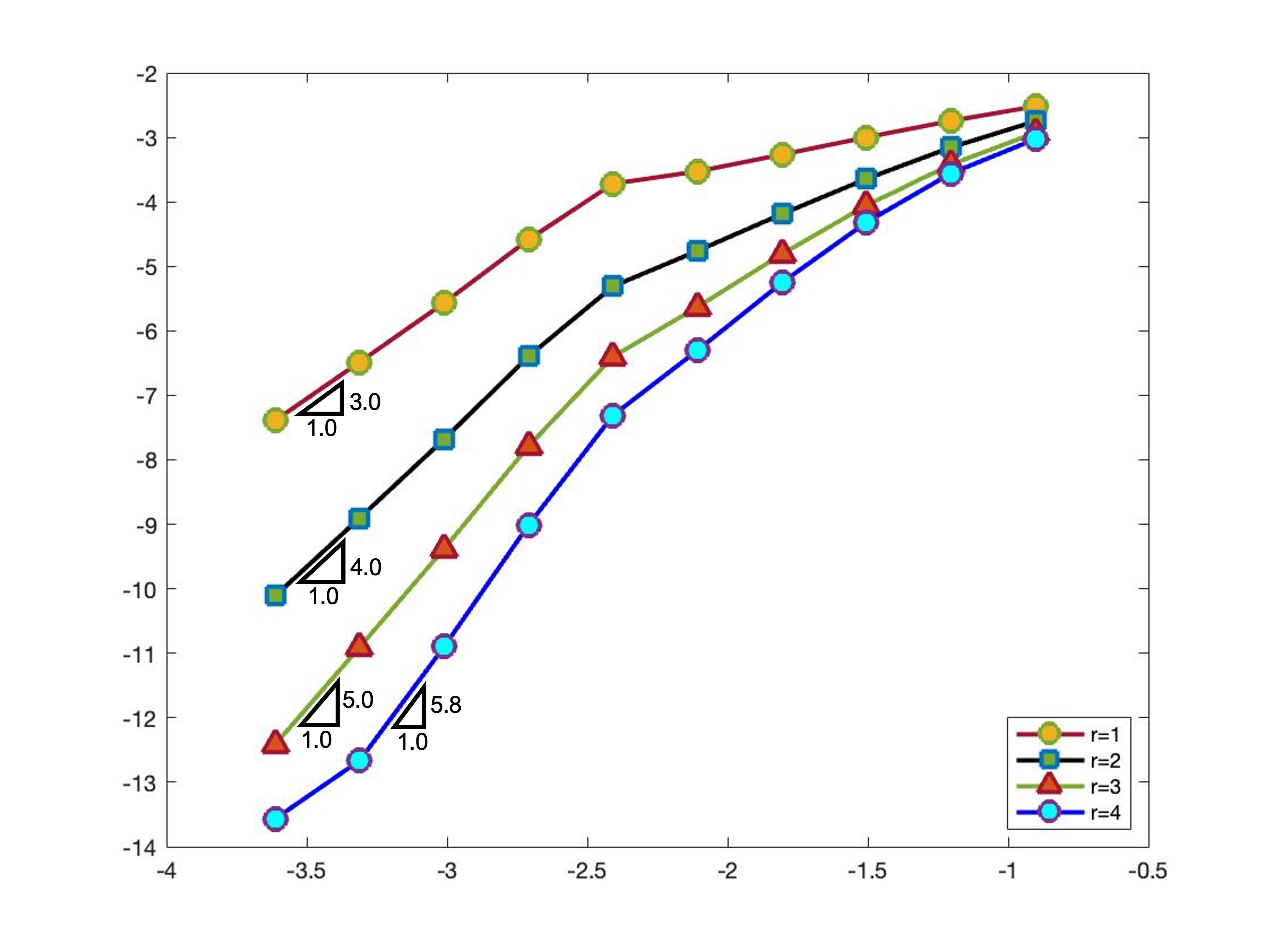} 
\label{fig:-exmp1b}
} 
\subfigure[$k=1000$]{\includegraphics[clip=true, trim=240 130 240 140, width=0.31\textwidth]{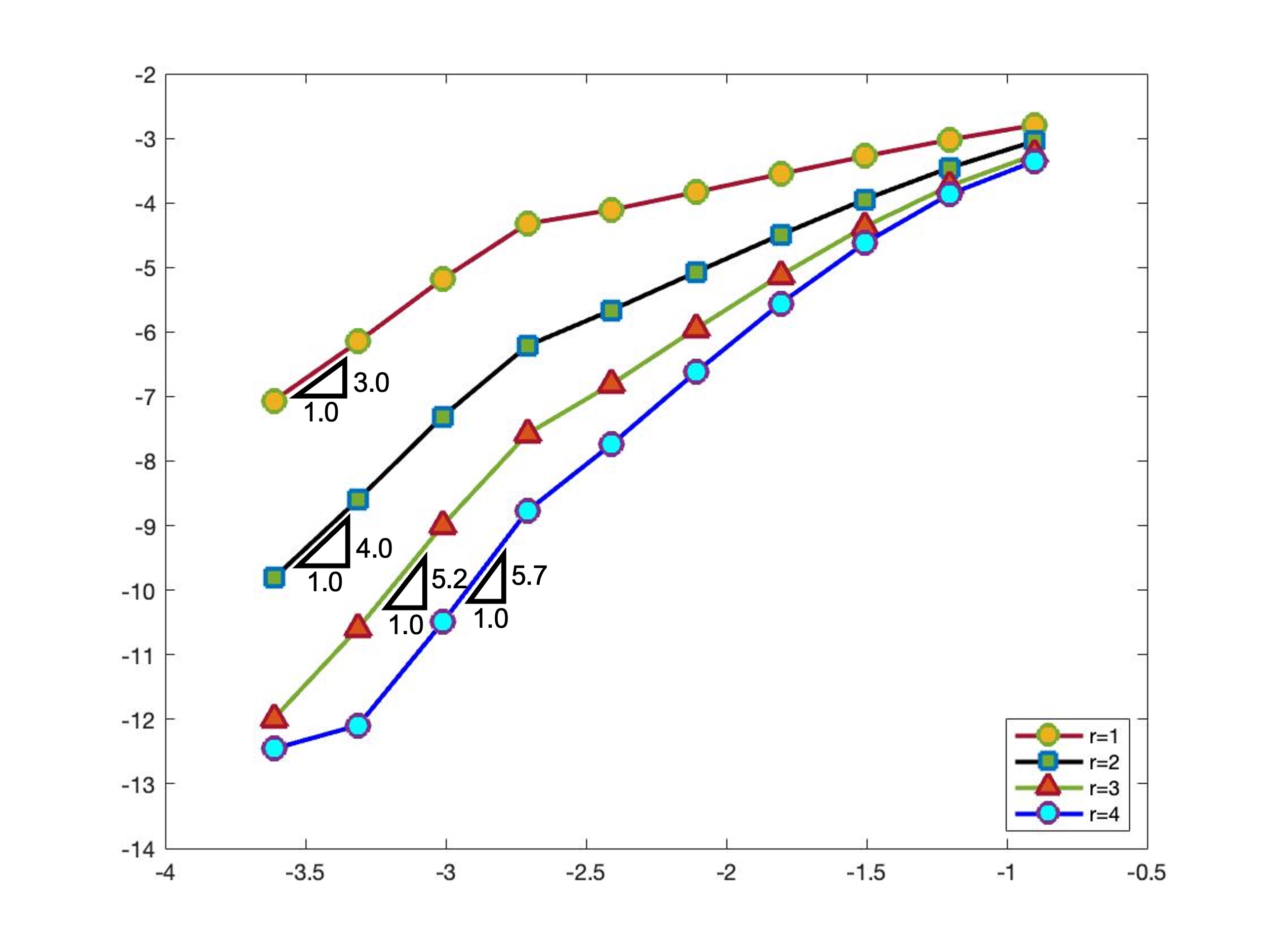} 
\label{fig:-exmp1c}
} 
\end{center}
\caption{The plot of $\log_{10} |I_{k,n}-I_k|$ against $\log_{10} \left(1/n\right)$ for integral in \cref{Eg1}}
\label{fig:exmp1}
\end{figure}


\end{exmp}

\begin{exmp} \label{Eg2} For studying the convergence behaviour of the numerical integration scheme on a second example without any singularity, that is, $w(x) = 1$, we consider the integral 
\begin{align*}
I_k = \int_{0}^{1}\sin(\cos t)\sin t \,e^{ik(\cos t)}\,dt
\end{align*}
The results in \cref{fig:exmp2} show that the order of convergence $r+2$ does not degrade with increasing oscillation frequency.

\begin{figure}[b]
\begin{center}
\subfigure[$k=10$ ]{\includegraphics[clip=true, trim=240 130 240 140, width=0.311\textwidth]{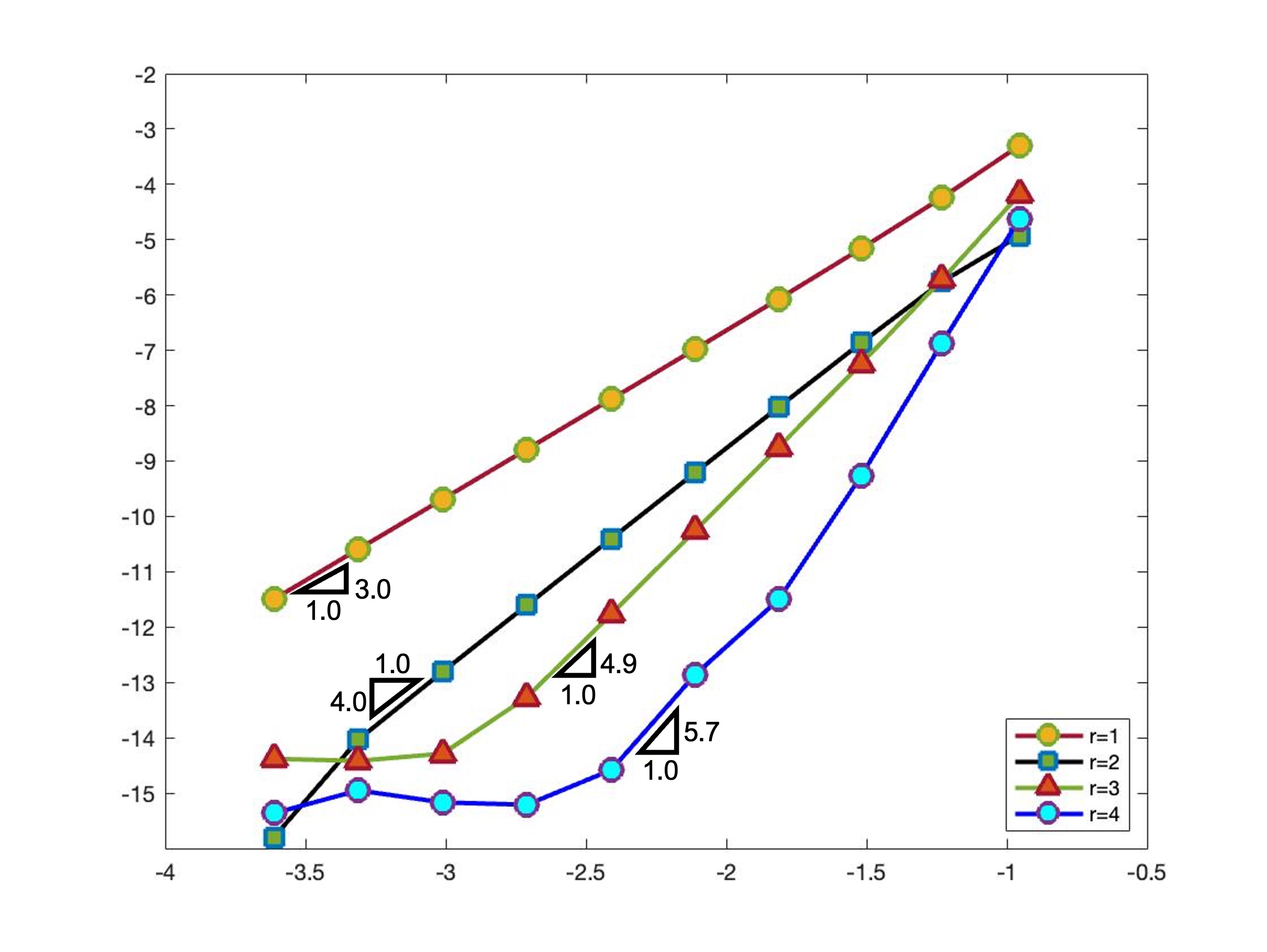} 
\label{fig:exmp2a}
} 
\subfigure[$k=100$]{\includegraphics[clip=true, trim=240 130 240 140, width=0.311\textwidth]{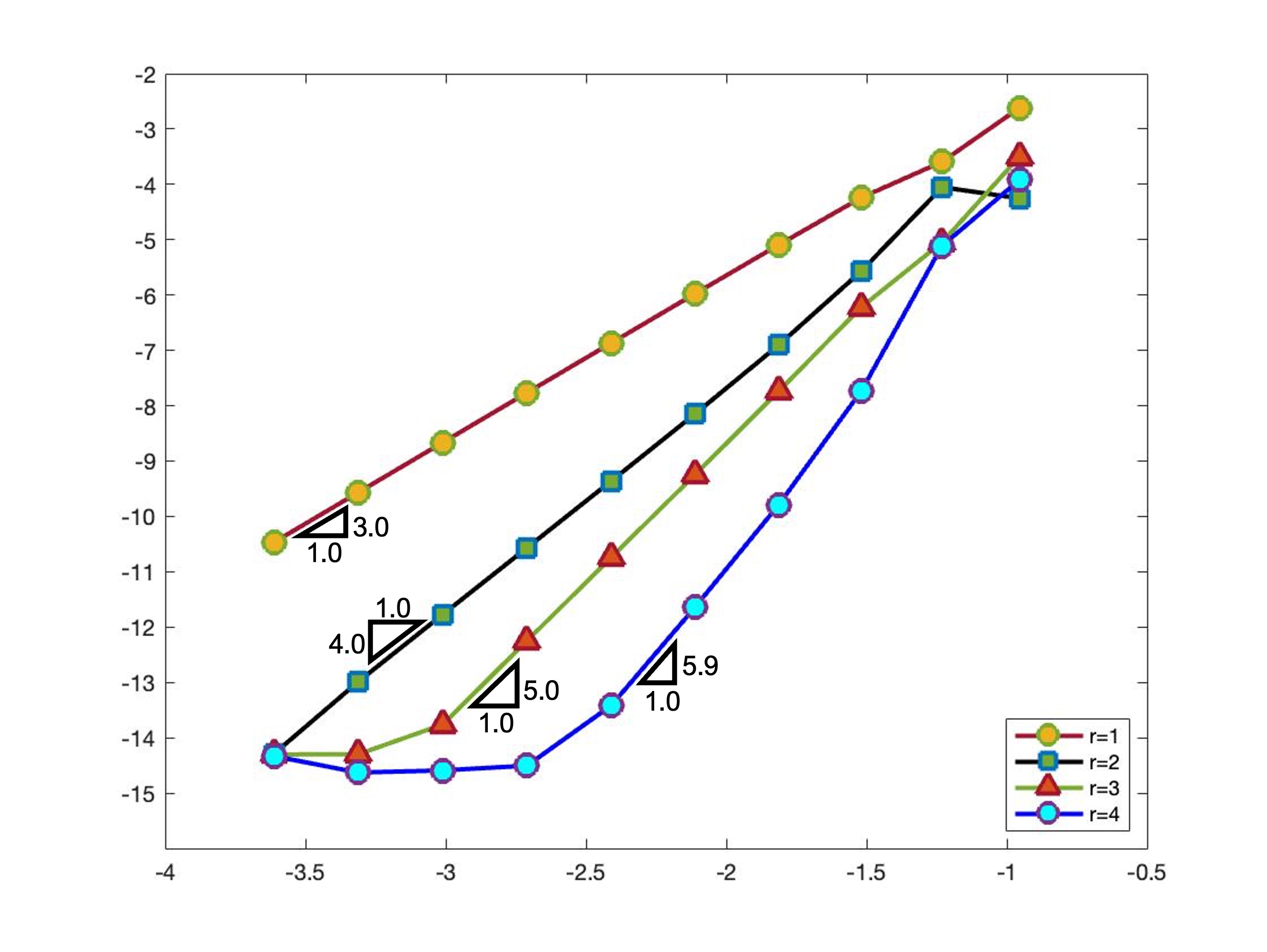} 
\label{fig:-exmp2b}
} 
\subfigure[$k=1000$]{\includegraphics[clip=true, trim=240 130 240 140, width=0.311\textwidth]{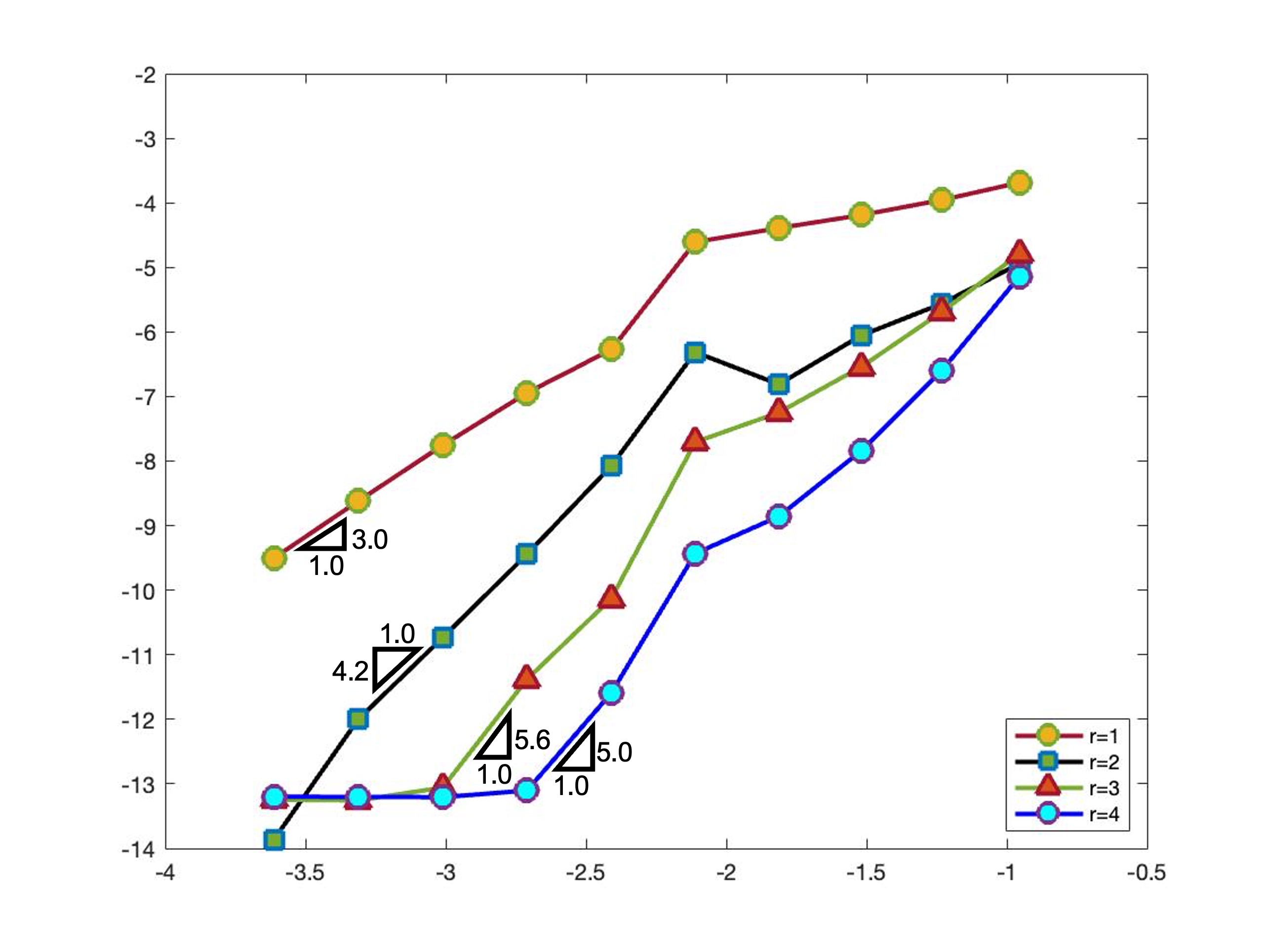} 
\label{fig:-exmp2b}
} 
\end{center}
\caption{The plot of $\log_{10} |I_{k,n}-I_k|$ against $\log_{10}(1/ n)$ for integral in \cref{Eg2}}
\label{fig:exmp2}
\end{figure}

\end{exmp}

\begin{exmp} \label{Eg3} For the next example, we consider the integral 
\begin{align*}
I_k^p = \int_{0}^{1}e^{ik t^p}\,dt \quad  = \frac{1}{p}  \int_{0}^{1} x^{-(p-1)/p} e^{ik x}\,dx 
\end{align*}
that has been studied in literature for measuring the performance of quadratures in dealing with end-point algebraic singularities. 
For instance, the results for this example is discussed in \cite{zaman2015new} where, for $k = 10000$ and $p = 2$,  a relative error of the order $10^{-7}$ is achieved with $n=100$. For $k = 3000$ and $p = 10$, a relative error of $10^{-12}$ is obtained with $n=100$. Approximations to these integrals using Filon-Clenshaw-Curtis have also been reported  in \cite{dominguez2013filon} for $p=2/3$ and $p=4/3$ under experiment 1.
Our approach, by construction, captures these integrals exactly upto the machine precision through analytic moment computations as seen in \cref{tab:exmp3} where we report the errors $|I^p_{k,2}-I^p_k|$ for $k=10^m,\ m = 3, 4, 5, 6, 7$ and $p = 2/3, 4/3, 2, 10$. The results shown are for $n=2$ and $r=0$. 

\begin{table} [!t] 
	\centering
	\begin{tabularx}{0.8\textwidth}{ >{\setlength\hsize{0.6\hsize}\centering}X | >{\setlength\hsize{1.1\hsize}\centering}X | >{\setlength\hsize{1.1\hsize}\centering}X | >{\setlength\hsize{1.1\hsize}\centering}X | >{\setlength\hsize{1.1\hsize}\centering}X  }
		\hline
		$k$ & $p=2/3$ & $p=4/3$ &  $p=2$  & $p=10$ \tabularnewline
		\hline\hline
		$10^3$ & $7.4325 \times 10^{-17} $ & $1.7110 \times 10^{-16}$ & $1.2337 \times 10^{-16} $ & $4.6653 \times 10^{-16}$
		\tabularnewline
		\hline
		$10^4$ & $0.0000 \times 10^{+00} $ & $2.7730 \times 10^{-16} $ & $9.7618 \times 10^{-17} $ & $5.8885 \times 10^{-16} $    
		\tabularnewline
		\hline
		$10^5$ & $2.2818 \times 10^{-17} $ & $2.2485 \times 10^{-16} $ & $1.5455 \times 10^{-16} $ & $5.5786 \times 10^{-16} $   
		\tabularnewline
		\hline
		$10^6$ & $0.0000 \times 10^{+00} $ & $2.9916 \times 10^{-16} $ & $0.0000 \times 10^{+00} $ & $8.1510 \times 10^{-16} $   
		\tabularnewline
		\hline
		$10^7$ & $0.0000 \times 10^{+00} $ & $6.9014 \times 10^{-16} $ & $1.3676 \times 10^{-16} $ & $4.4208 \times 10^{-16} $   
		\tabularnewline
		
		\hline
	\end{tabularx}
	\caption{The approximation errors $|I_{k,n}^p-I_k^p|$ for the integral in \cref{Eg3} with $n=2$ and $r=0$}
	\label{tab:exmp3}
\end{table}

 \end{exmp}

\begin{exmp} \label{Eg4} Next, we study the performance of our method in approximating the integral
\begin{align*}
I_k = \int_{0}^{\pi} e^{ik  \sqrt{(\cos t - 1)^2+ \sin^2 t}}\,dt,
\end{align*}
Using the change of variable $x = \sqrt{(\cos t - 1)^2+ \sin^2 t}$, we have 
\begin{align*}
I_k = \int_{0}^{2} \frac{2}{ \sqrt{4-x^2}} e^{ik  x}\,dx = \int_{0}^{2} \frac{1}{ \sqrt{2-x}} \frac{2}{ \sqrt{2+x}} e^{ik  x}\,dx,
\end{align*}
where we have an algebraic singularity of order $\beta = -1/2$ at the upper end point of the integration interval. Moreover, the function $f(x) = 1/\sqrt{2+x}$ is smooth on $[0,2]$. The results in \cref{fig:exmp4} show that the approximations convergence at the expected rate  $r + 3/2$.

\begin{figure}[b]
\begin{center}
\subfigure[$k=100$ ]{\includegraphics[clip=true, trim=240 130 240 140, width=0.311\textwidth]{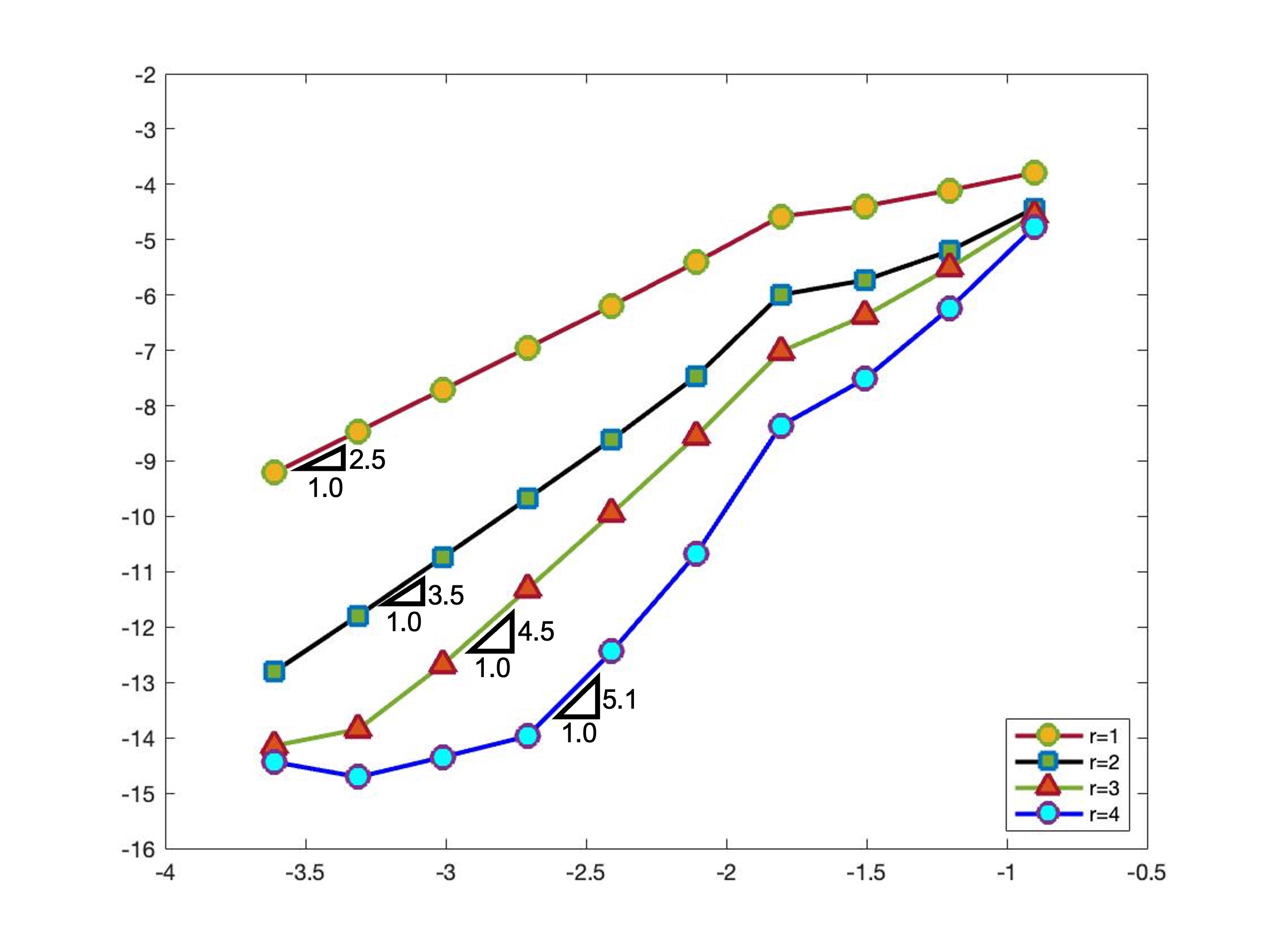} 
\label{fig:exmp4a}
} 
\subfigure[$k=500$]{\includegraphics[clip=true, trim=240 130 240 140, width=0.311\textwidth]{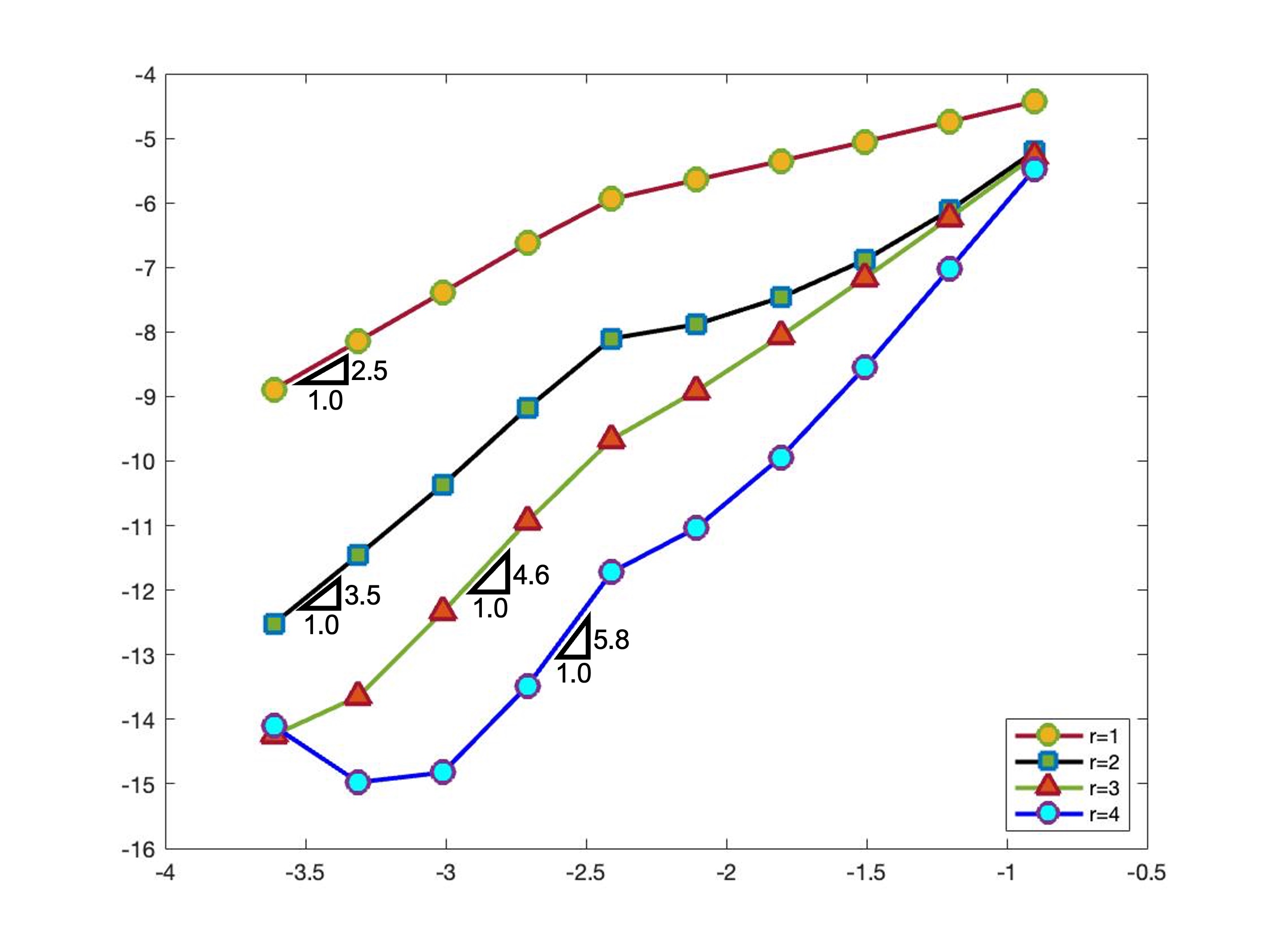} 
\label{fig:-exmp4b}
} 
\subfigure[$k=1000$]{\includegraphics[clip=true, trim=240 130 240 140, width=0.311\textwidth]{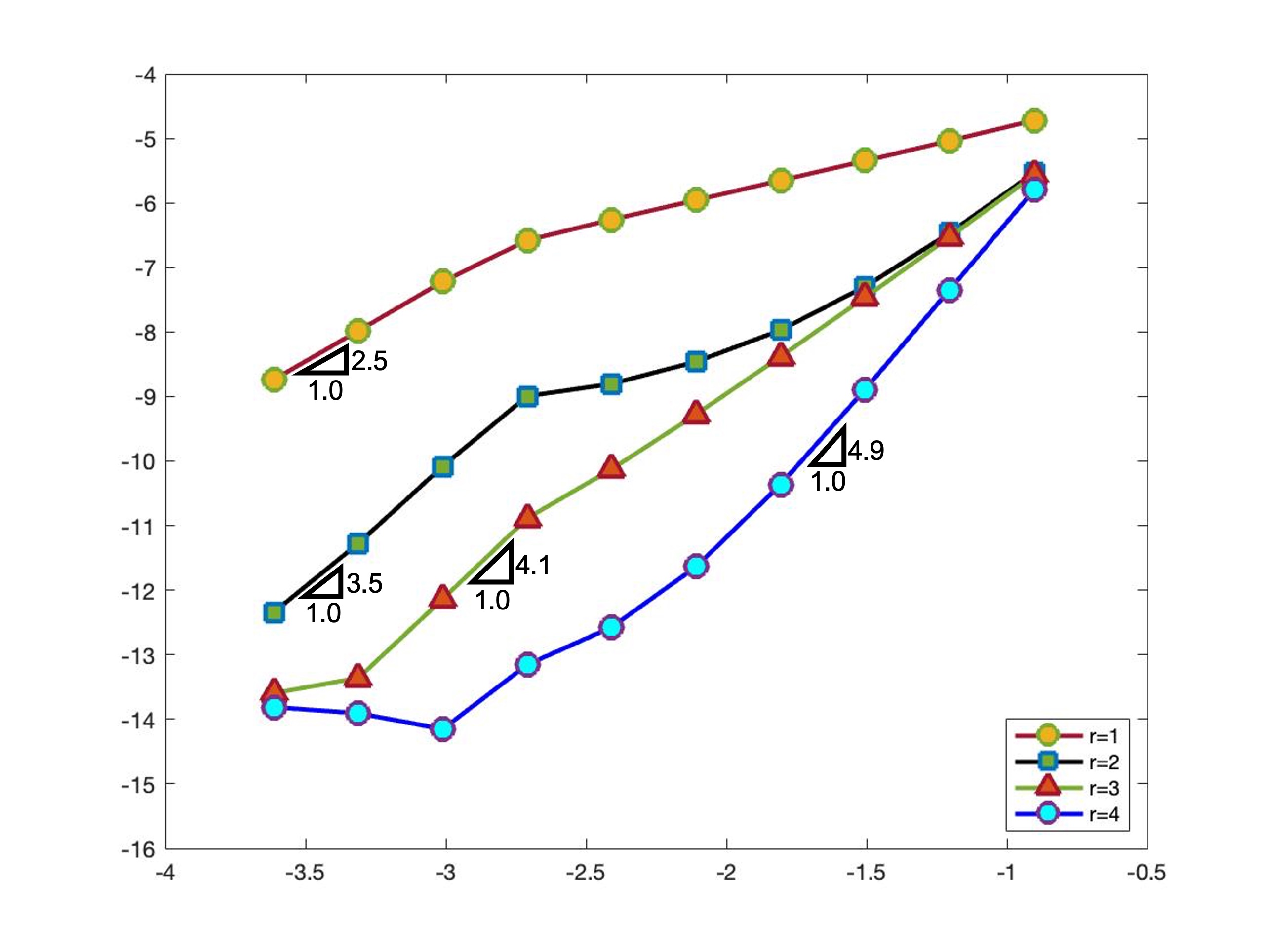} 
\label{fig:-exmp4b}
} 
\end{center}
\caption{The plot of $\log_{10} |I_{k,n}-I_k|$ against $\log_{10}(1/ n)$ for integral in \cref{Eg4}}
\label{fig:exmp4}
\end{figure}

\end{exmp}

\begin{exmp} \label{Eg5}
For the following integral
\begin{align*}
I_k = \int_{0}^{1} x^\beta \left|x-\frac{1}{2}\right|  e^{ik x}\,dx
\end{align*}
where $f(x)=|x-1/2|$ is only piecewise smooth, high-order approximation can still be obtained by breaking the integral appropriately as seen in the results reported in \cref{fig:exmp5} where numerical integrals exhibit the convergence rate $r+2+\beta$.

\begin{figure}[t]
\begin{center}
\subfigure[$\beta = -1/2$ ]{\includegraphics[clip=true, trim=240 105 220 150, width=0.311\textwidth]{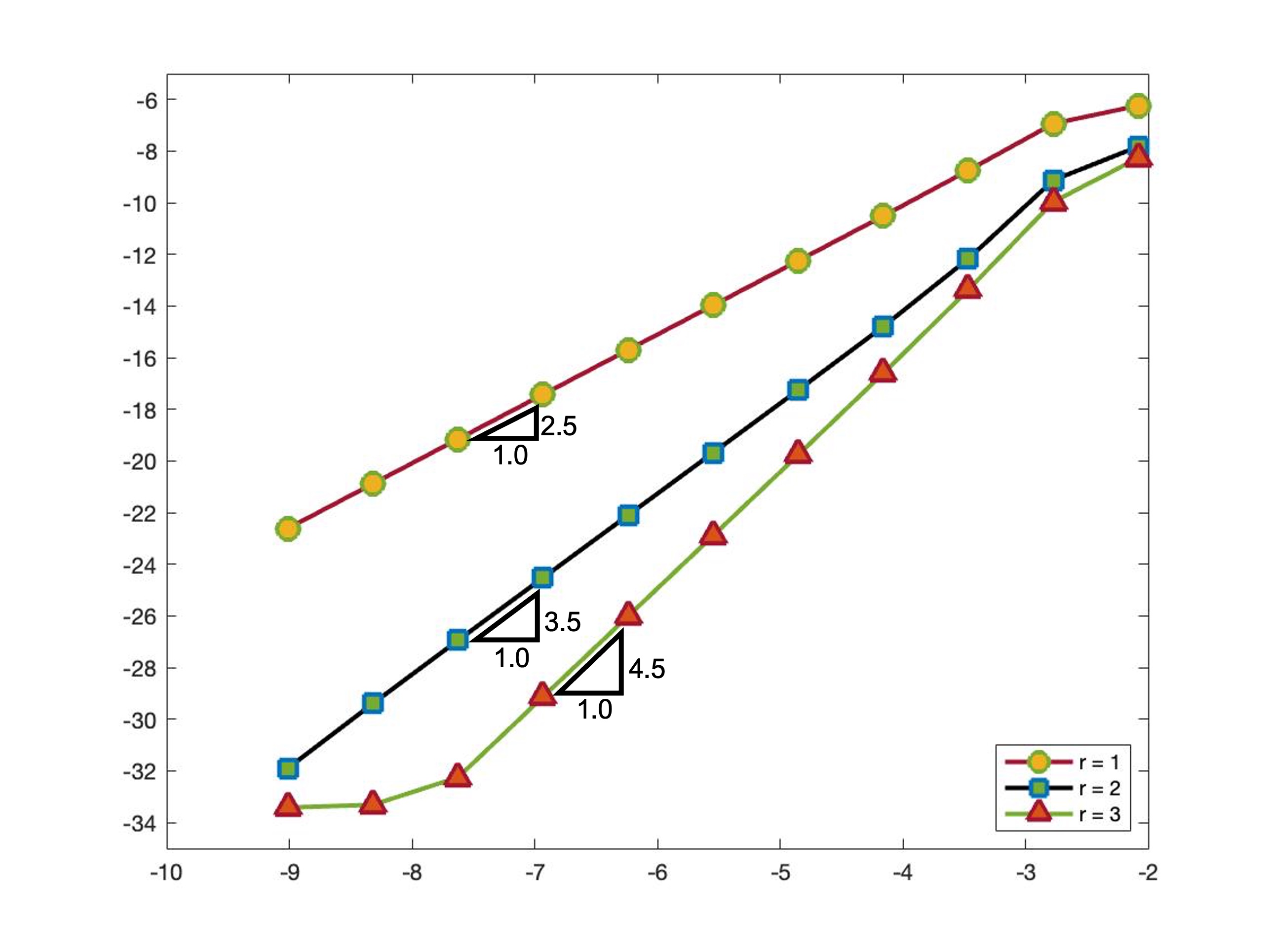} 
\label{fig:exmp5a}
} 
\subfigure[$\beta = -1/4$]{\includegraphics[clip=true, trim=240 105 220 150, width=0.311\textwidth]{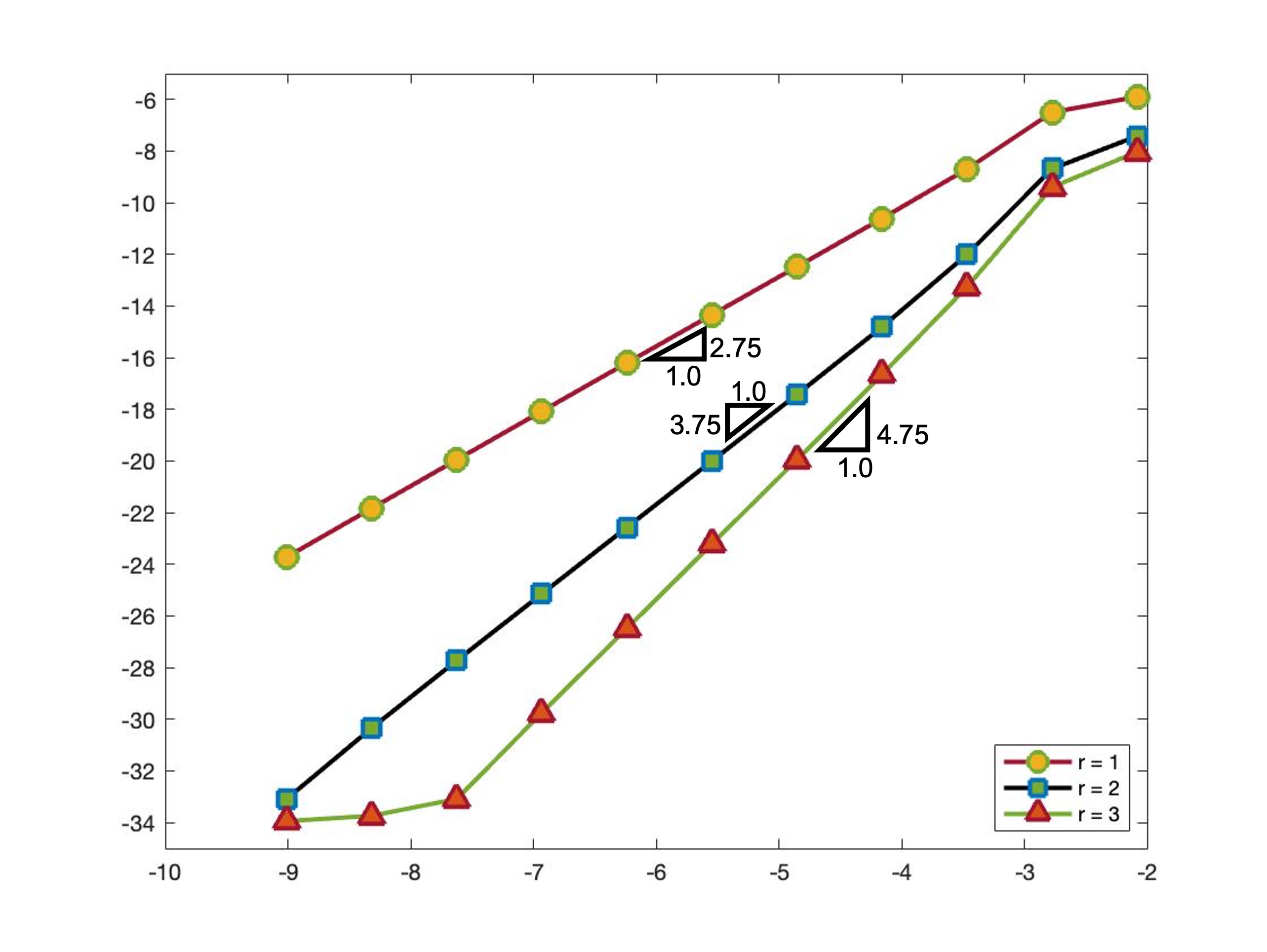} 
\label{fig:-exmp5b}
} 
\subfigure[$\beta = -2/3$]{\includegraphics[clip=true, trim=240 105 220 150, width=0.311\textwidth]{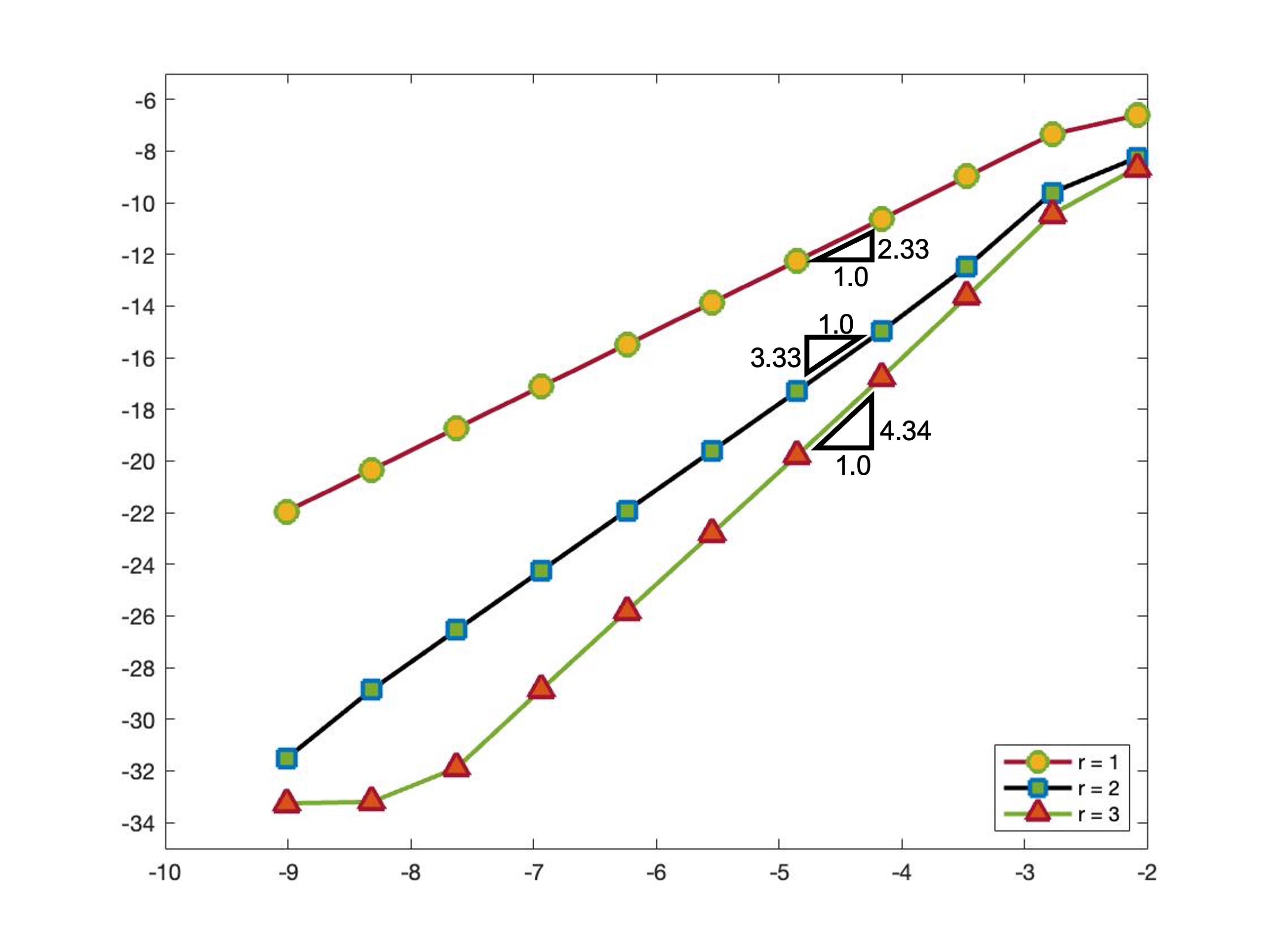} 
\label{fig:-exmp5c}
} 
\end{center}
\caption{The plot of $\log_{10} |I_{k,n}-I_k|$ against $\log_{10}(1/ n)$ for integral in \cref{Eg5}}
\label{fig:exmp5}
\end{figure}

\end{exmp}

\begin{exmp} \label{Eg6} We now consider an integral with a weight of the form $w(x) = (x-a)^\alpha (b-x)^\beta$. In particular, we approximate the following integral
\begin{align*}
I_k^1 = \int_{0}^{1}x^{-1/2}(1-x)^{-1/3}e^x e^{ik x}\,dx
\end{align*}
using the proposed scheme for $k=10$ and $k=100$. The convergence rates for $r=1,2,3,4$, shown in \cref{fig:exmp61}, are consistent with the theoretical rate $r+2-\max\{-\alpha,-\beta\}$.

\begin{figure}[t]
\begin{center}
\subfigure[$k = 10$ ]{\includegraphics[clip=true, trim=240 105 220 150, width=0.311\textwidth]{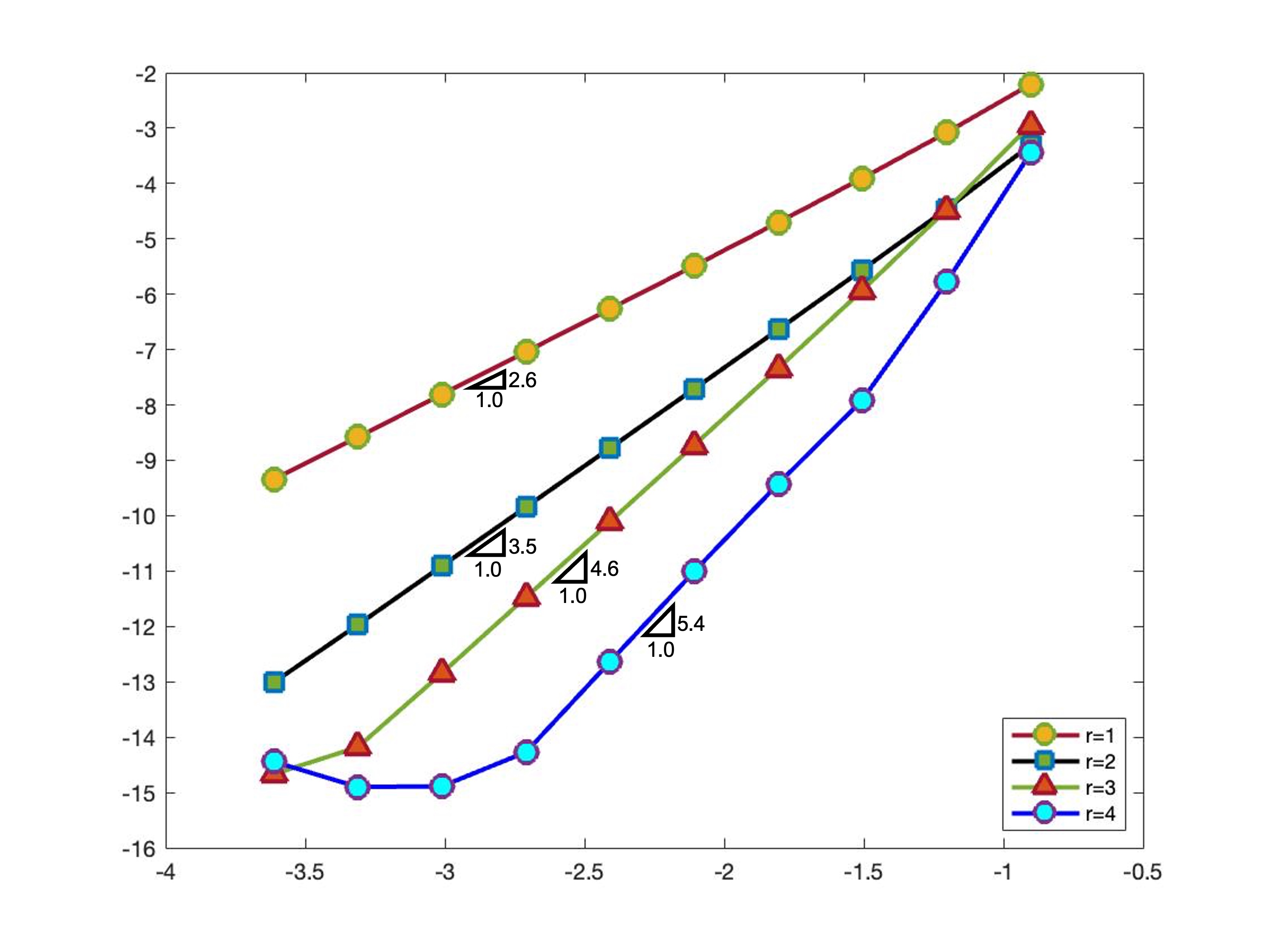} 
\label{fig:exmp5a}
} 
\subfigure[$k = 100$]{\includegraphics[clip=true, trim=240 105 220 150, width=0.311\textwidth]{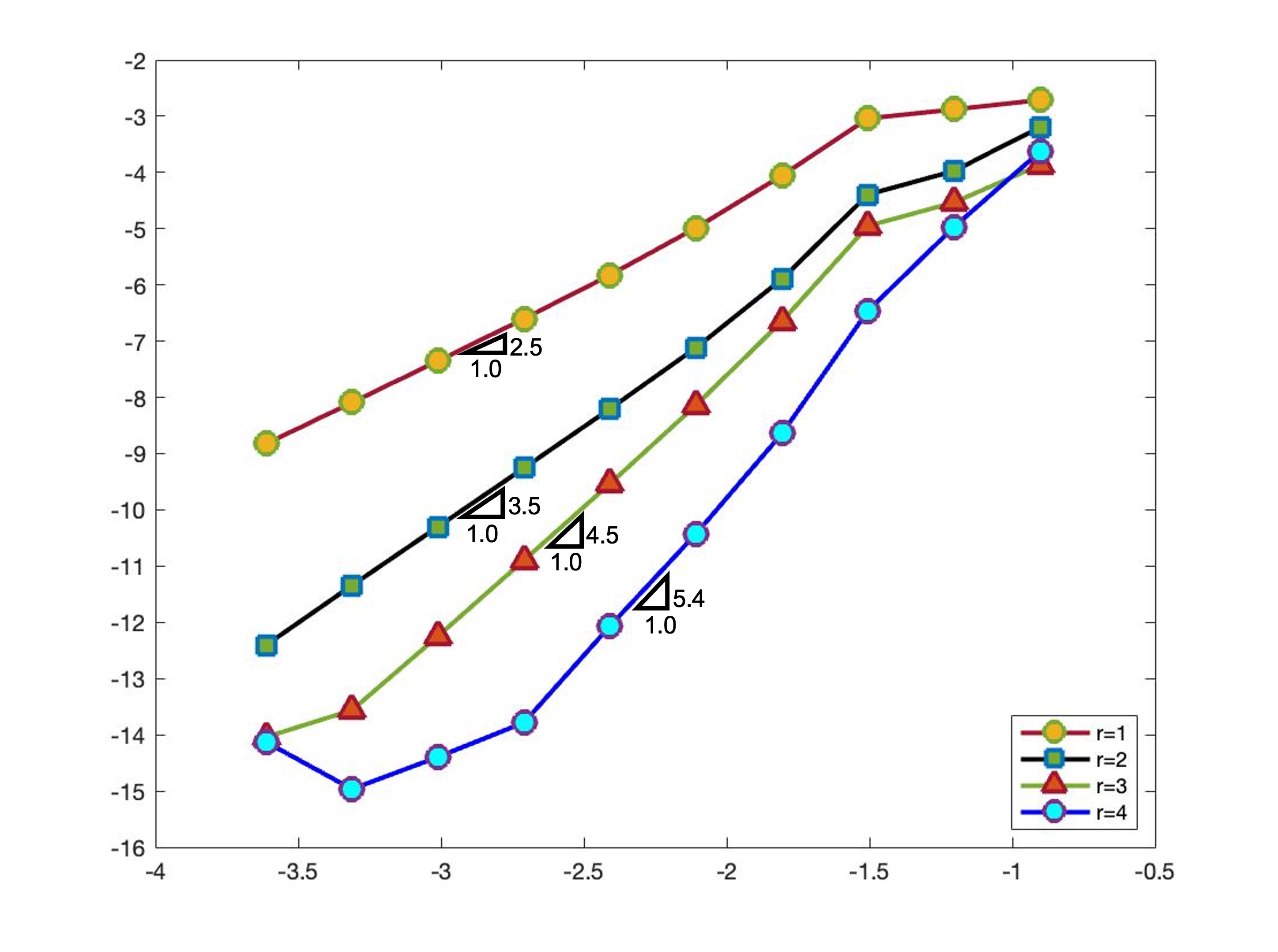} 
\label{fig:-exmp5b}
} 
\subfigure[$k = 500$]{\includegraphics[clip=true, trim=240 105 220 150, width=0.311\textwidth]{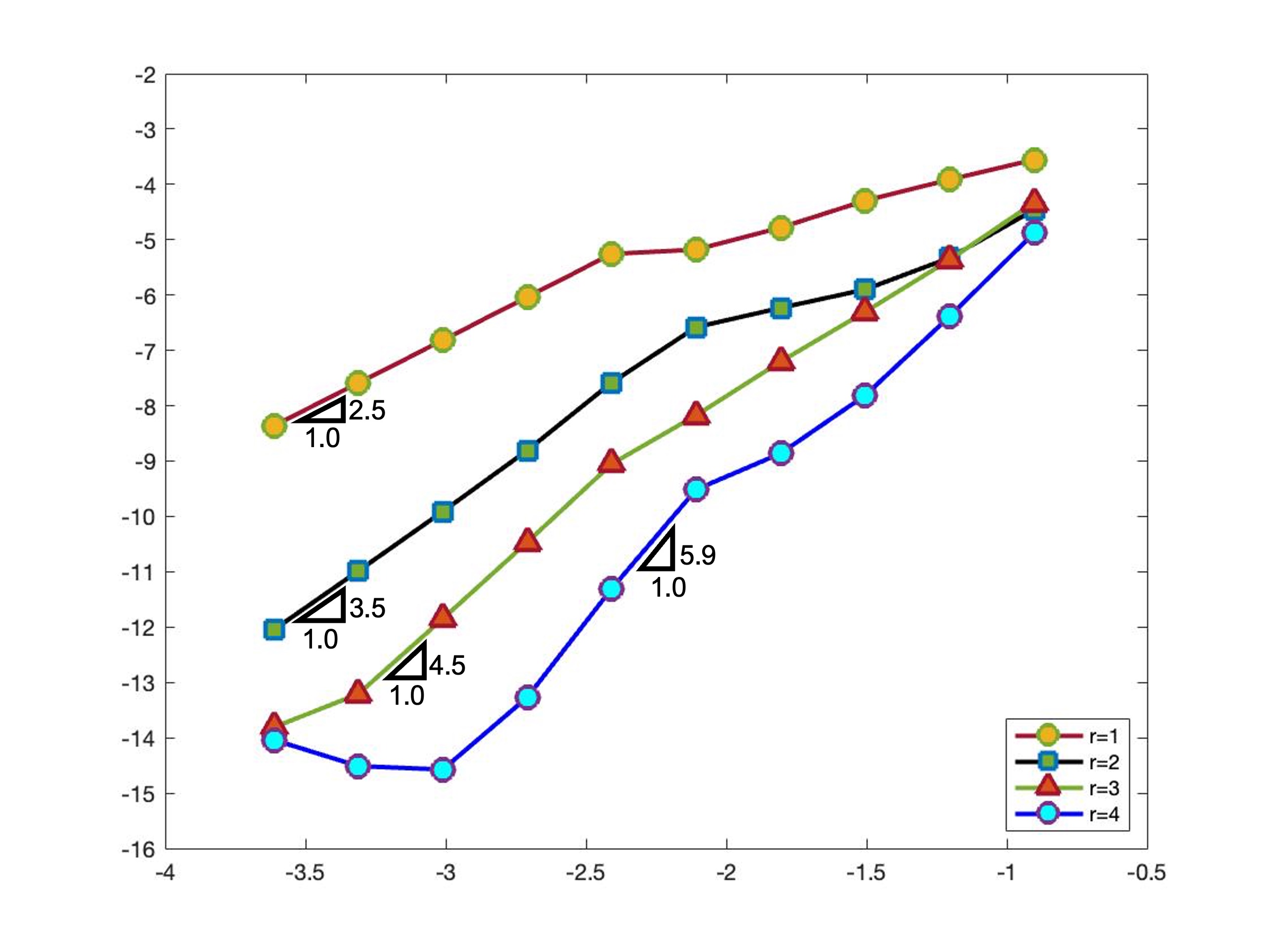} 
\label{fig:-exmp5c}
} 
\end{center}
\caption{The plot of $\log_{10} |I_{k,n}^1-I_k^1|$ against $\log_{10}(1/ n)$ for integral in \cref{Eg6}}
\label{fig:exmp61}
\end{figure}


The the results for a second integral of this type
\begin{align*}
I_k^2= \int_{2}^{3}(x-2)^{-1/4}(3-x)^{-2/3}\sin(x) e^{ik x}\,dx,
\end{align*}
presented in \cref{fig:exmp62}, confirm the convergence rates.

\begin{figure}[t]
\begin{center}
\subfigure[$k = 10$ ]{\includegraphics[clip=true, trim=240 105 220 150, width=0.311\textwidth]{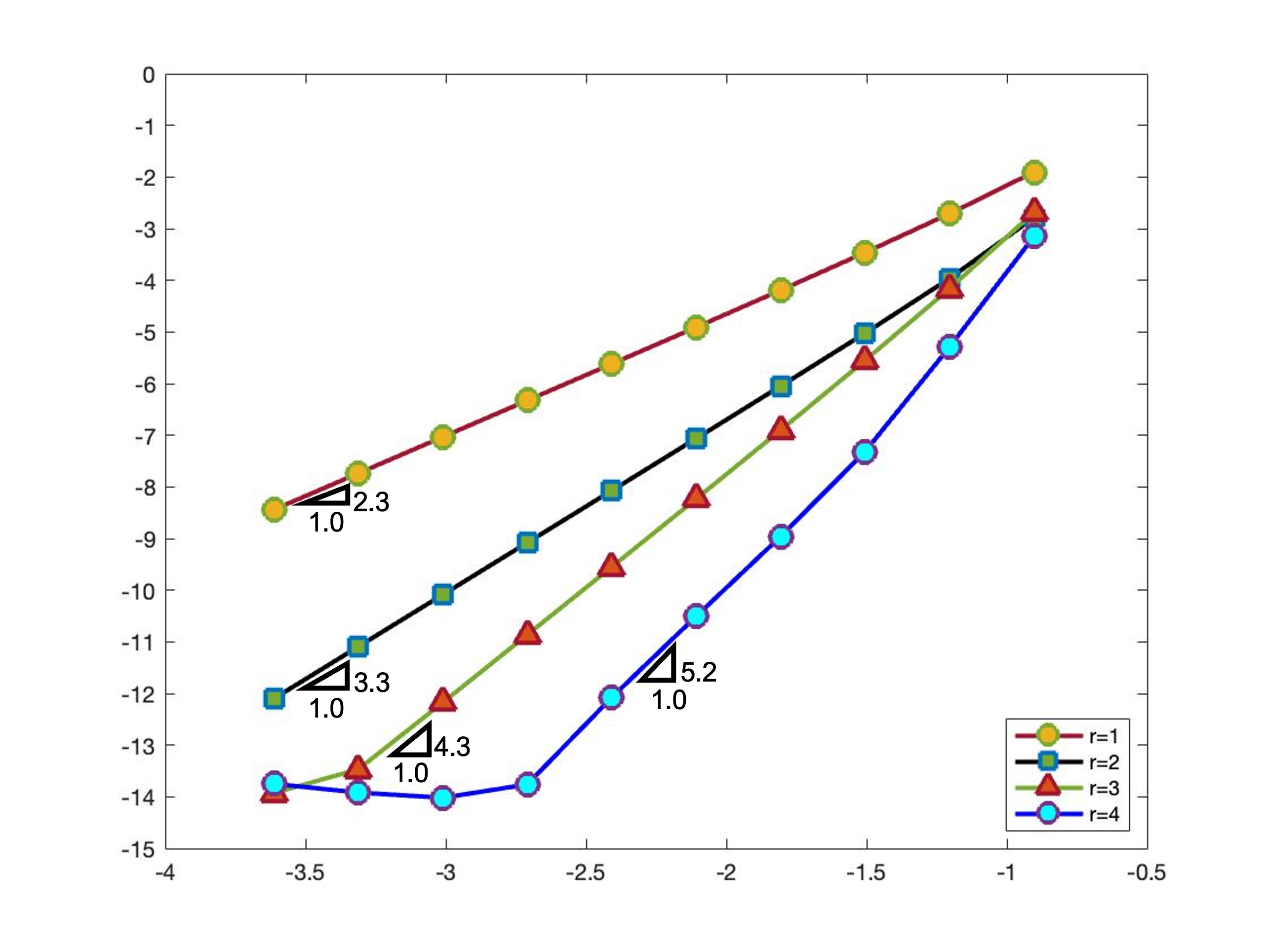} 
\label{fig:exmp5a}
} 
\subfigure[$k = 100$]{\includegraphics[clip=true, trim=240 105 220 150, width=0.311\textwidth]{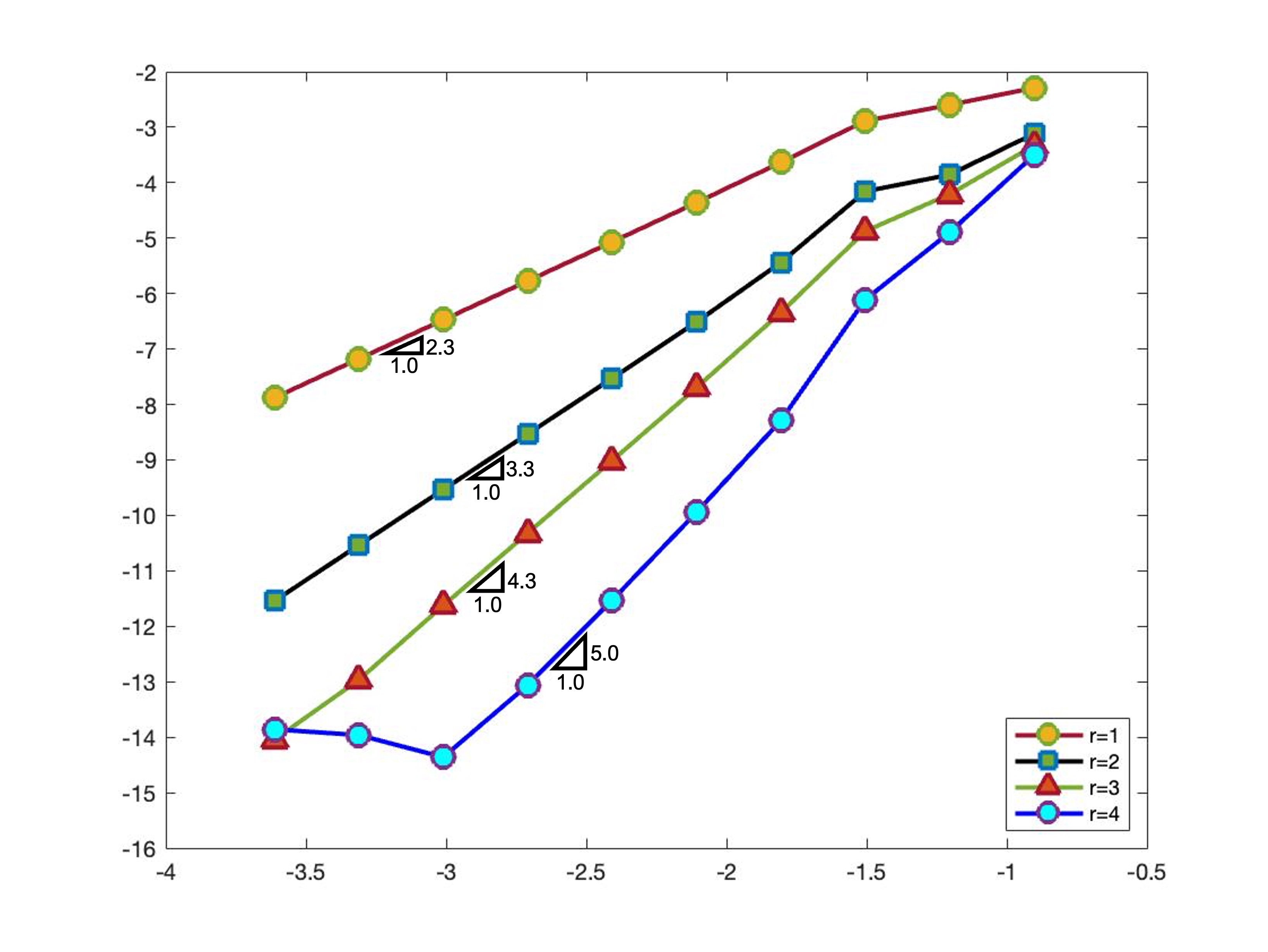} 
\label{fig:-exmp5b}
} 
\subfigure[$k = 500$]{\includegraphics[clip=true, trim=240 105 220 150, width=0.311\textwidth]{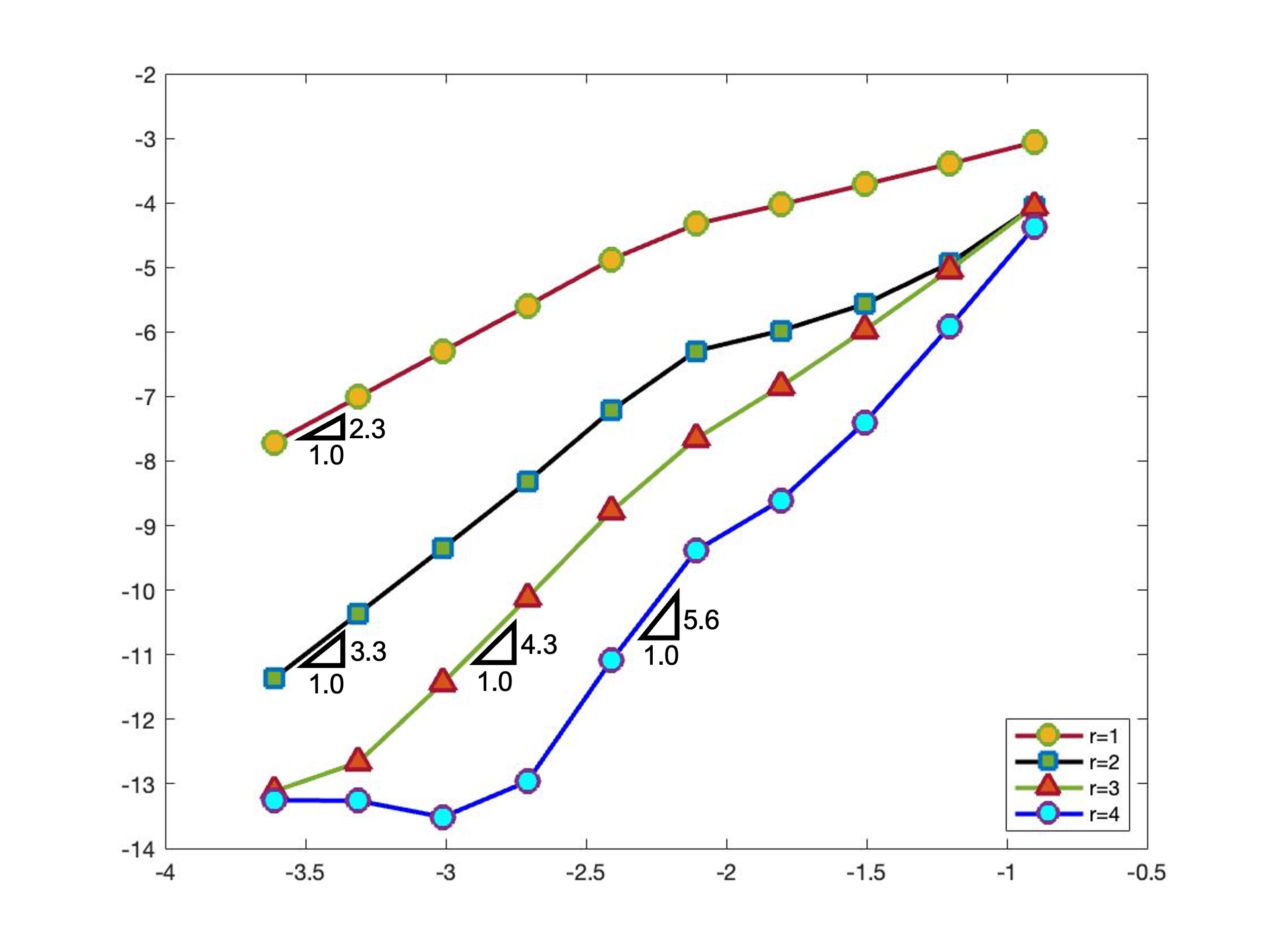} 
\label{fig:-exmp5c}
} 
\end{center}
\caption{The plot of $\log_{10} |I_{k,n}^2-I_k^2|$ against $\log_{10}(1/ n)$ for integral in \cref{Eg6}}
\label{fig:exmp62}
\end{figure}


\end{exmp}

\begin{exmp} \label{Eg10} Finally, we study the performance of the proposed method in approximating the following integral 
\begin{align*}
I_k = \int_{0}^{\pi/2} \log\left( \sqrt{(\cos t - 1)^2+ \sin^2 t} \right)e^{ik  \sqrt{(\cos t - 1)^2+ \sin^2 t}}\,dt,
\end{align*}
that has a logarithmic singularity at $t=0$. As in \cref{Eg4}, using the change of variable $x = \sqrt{(\cos t - 1)^2+ \sin^2 t}$, yields a linear oscillator
\begin{align*}
I_k = \int_{0}^{\sqrt{2}} \log(x) \frac{2}{ \sqrt{4-x^2}} e^{ik  x}\,dx ,
\end{align*}
where $w(x) = \log(x)$ and $f(x) = 2/\sqrt{4-x^2}$ is smooth on $[0,\sqrt{2}]$.

\begin{figure}[t]
\begin{center}
\subfigure[$k = 100$ ]{\includegraphics[clip=true, trim=240 105 220 150, width=0.311\textwidth]{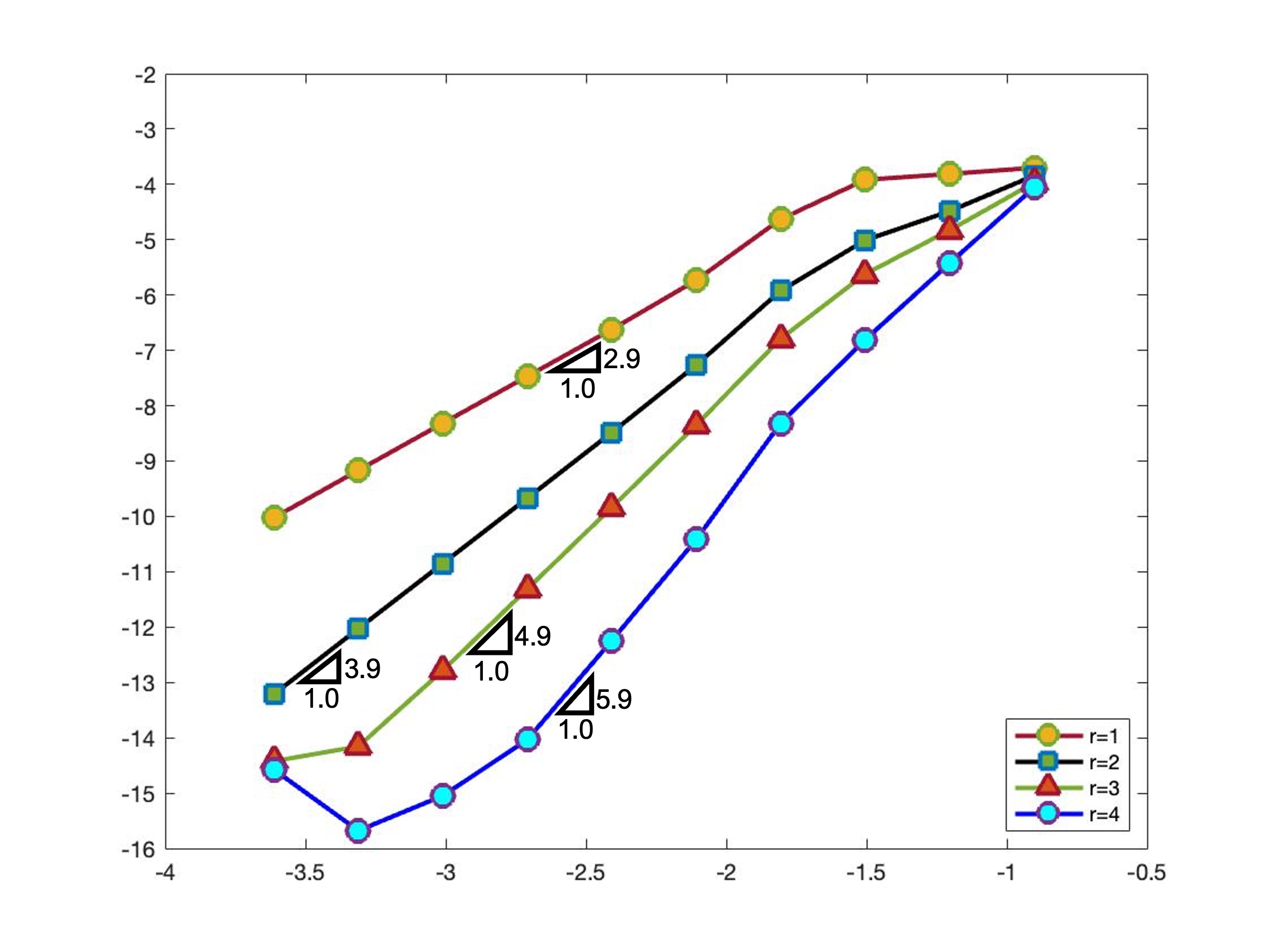} 
\label{fig:exmp5a}
} 
\subfigure[$k = 500$]{\includegraphics[clip=true, trim=240 105 220 150, width=0.311\textwidth]{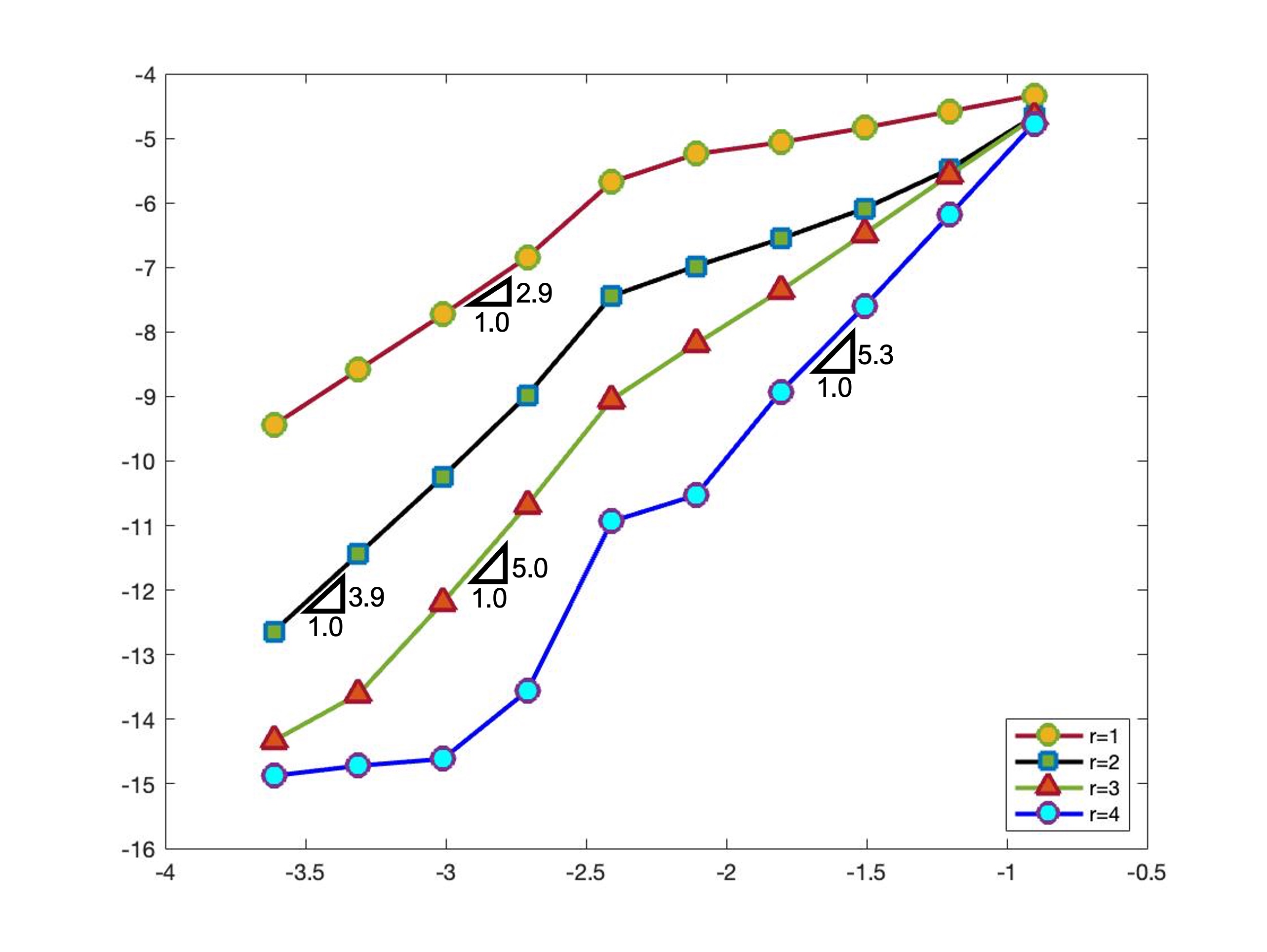} 
\label{fig:-exmp5b}
} 
\subfigure[$k = 1000$]{\includegraphics[clip=true, trim=240 105 220 150, width=0.311\textwidth]{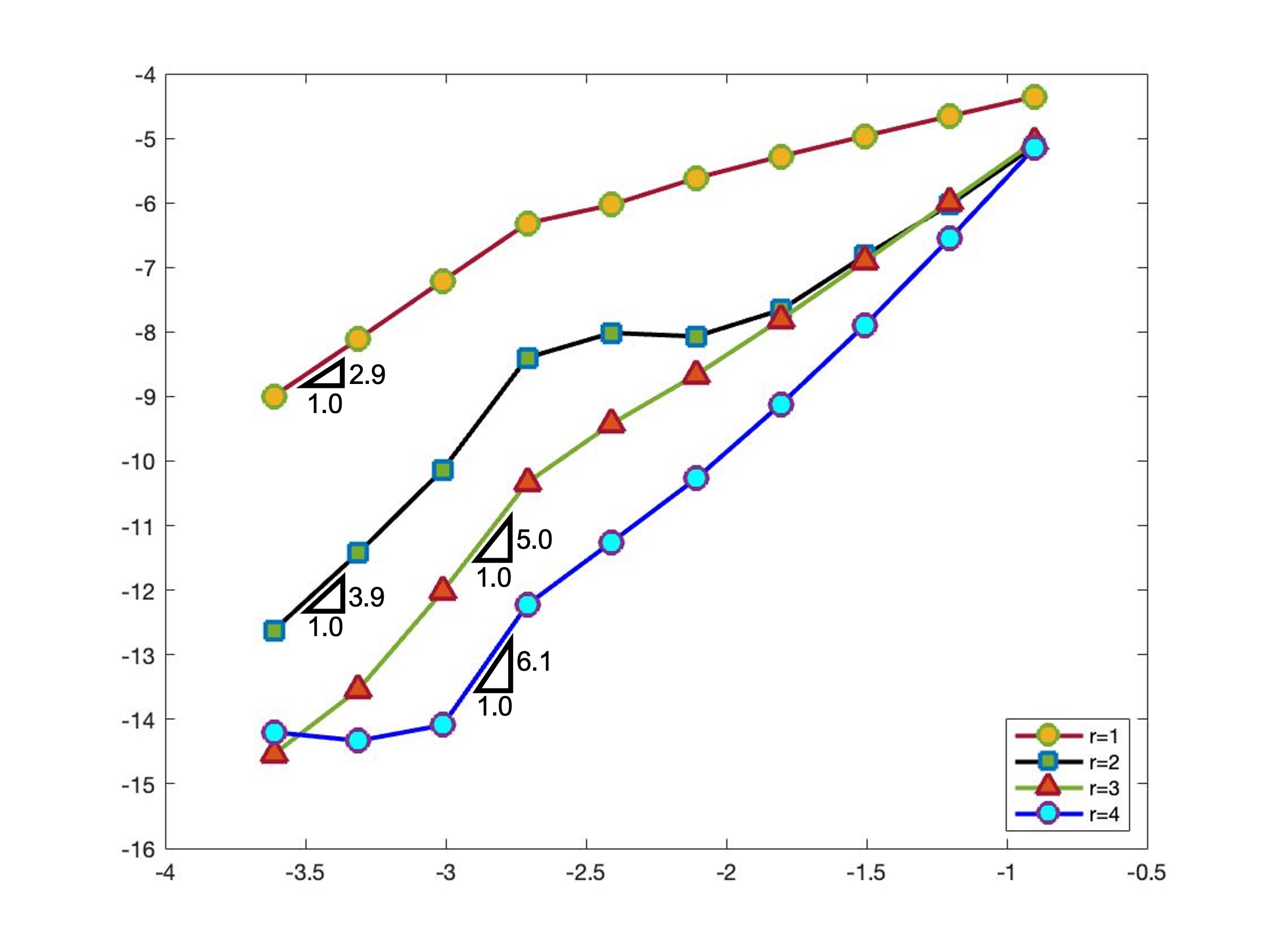} 
\label{fig:-exmp5c}
} 
\end{center}
\caption{The plot of $\log_{10} |I_{k,n}-I_k|$ against $\log_{10}(1/n)$ for integral in \cref{Eg10}}
\label{fig:exmp5}
\end{figure}

%
%

\end{exmp}

%
%
%

\section{Conclusion}

In this paper, we introduce a new Filon-type quadrature for numerical of highly oscillatory integrals of the form  \cref{WeightGenProb} that can be employed to approximate more general oscillatory integrals of the form \cref{genprob} with algebraic singularities and stationary points. The proposed method utilises an interpolating trigonometric polynomial to approximate the smooth and non-oscillatory part of the integrand. Normally, when a function is not periodic, such trigonometric polynomial approximations suffer from the Gibb's phenomenon. In order that the unwanted Gibb's oscillations do not adversely impact the accuracy of numerical integration, the method constructs a periodic extension of the underlying function before its approximation. The main advantage of using  trigonometric interpolation lies in obtaining a relatively simpler moment problem that can be solved analytically for many important classes of integrands. This approach can easily handle one and two sided integrable algebraic singularities as well as logarithmic singularities, The accompanying numerical analysis of the method reveals that the approximations converge rapidly to the exact integral. The scheme is tested on a variety of highly oscillatory integrals and these numerical experiments confirm that the approximations converge at the rates predicted by the theory.  

\bibliography{mybibliography.bib}
\bibliographystyle{elsarticle-num} 

\end{document}